\documentclass[a4paper]{amsart}
\usepackage[T1]{fontenc}
\usepackage[latin1]{inputenc}
\usepackage{amssymb}
\usepackage{amsfonts}
\usepackage{amsxtra}
\usepackage{ae}
\usepackage[all]{xy}
\usepackage{enumerate}
\usepackage[lite]{amsrefs}

\newcommand*{\KK}{\mathrm{KK}}

\newcommand*{\K}{\mathrm{K}}
\newcommand*{\Ktop}{\mathrm{K}^{\mathrm{top}}}

\newcommand*{\HH}{\mathrm{HH}}
\newcommand*{\HC}{\mathrm{HC}}
\newcommand*{\HP}{\mathrm{HP}}
\newcommand*{\CHH}{\underline{\mathrm{HH}}}
\newcommand*{\CHo}{\underline{\mathrm{H}}}

\newcommand*{\fin}{{\mathrm{f}}}
\newcommand*{\dR}{\mathrm{dR}}
\newcommand*{\inhomogen}{\mathrm{inh}}
\newcommand*{\homogen}{\mathrm{h}}
\newcommand*{\Ad}{{\mathrm{Ad}}}

\newcommand*{\abs}[1]{\lvert#1\rvert}

\newcommand*{\conj}[1]{\overline{#1}}
\newcommand*{\cl}[1]{\overline{#1}}
\newcommand*{\Cont}[1]{\complement #1}
\newcommand*{\inv}{^\times}

\newcommand*{\blank}{{\llcorner\!\!\lrcorner}}

\newcommand*{\GF}{{K}}

\newcommand*{\MCR}{{\mathcal{O}}}

\newcommand*{\Adel}{{\mathbb{A}}}
\newcommand*{\ICL}{{\mathcal{C}}}
\newcommand*{\Roots}{{\mu_\GF}}

\newcommand*{\Places}{\mathcal{P}}

\newcommand*{\Mod}{\mathbf{M}}
\newcommand*{\Der}{\mathbf{D}}

\newcommand*{\ID}{\mathrm{id}}

\newcommand*{\Sch}{\mathcal{S}}
\newcommand*{\Twist}{\mathcal{T}}
\newcommand*{\CCINF}{\mathcal{D}}

\newcommand*{\Fourier}{\mathfrak{F}}

\newcommand*{\hot}{\mathbin{\hat{\otimes}}}
\newcommand*{\Lhot}{\mathbin{\mathbb{L}{\hat{\otimes}}}}

\newcommand*{\cross}{\mathbin{\ltimes}}
\newcommand*{\defeq}{\mathrel{:=}}
\newcommand*{\into}{\rightarrowtail}
\newcommand*{\prto}{\twoheadrightarrow}
\newcommand*{\congto}{\overset{\cong}\to}

\DeclareMathOperator{\supp}{supp}

\DeclareMathOperator{\Ker}{Ker}

\DeclareMathOperator{\Hom}{Hom}

\DeclareMathOperator{\Rep}{Rep}

\DeclareMathOperator{\cInd}{c-Ind}
\DeclareMathOperator{\Res}{Res}

\newcommand*{\op}{\mathrm{op}}

\newcommand*{\C}{{\mathbb{C}}}
\newcommand*{\F}{{\mathbb{F}}}
\newcommand*{\R}{{\mathbb{R}}}
\newcommand*{\Z}{{\mathbb{Z}}}
\newcommand*{\Ztwo}{{\mathbb{Z}/2}}
\newcommand*{\N}{{\mathbb{N}}}
\newcommand*{\Q}{{\mathbb{Q}}}

\newcommand*{\brd}{-\hspace{0pt}}
\newcommand*{\nbd}{\nobreakdash-\hspace{0pt}}

\theoremstyle{plain}
\newtheorem{theorem}{Theorem}[section]
\newtheorem{proposition}[theorem]{Proposition}
\newtheorem{lemma}[theorem]{Lemma}

\theoremstyle{definition}

\theoremstyle{remark}

\begin{document}

\title[Cyclic homology of crossed products]{The cyclic homology and
  $\mathrm{K}$\nobreakdash-theory of certain adelic crossed products}
\author{Ralf Meyer}
\address{Mathematisches Institut\\
         Westfälische Wilhelms-Universität Münster\\
         Einsteinstr.\ 62\\
         48149 Münster\\
         Germany
}
\email{rameyer@math.uni-muenster.de}


\thanks{This research was supported by the EU-Network \emph{Quantum
    Spaces and Noncommutative Geometry} (Contract HPRN-CT-2002-00280)
  and the \emph{Deutsche Forschungsgemeinschaft} (SFB 478).}

\begin{abstract}
  The multiplicative group of a global field acts on its adele ring by
  multiplication.  We consider the crossed product algebra of the
  resulting action on the space of Schwartz functions on the adele
  ring and compute its Hochschild, cyclic and periodic cyclic
  homology.  We also compute the topological $\K$\nbd{}theory of the
  $C^*$\nbd{}algebra crossed product.
\end{abstract}
\maketitle

\section{Introduction}
\label{sec:intro}

Let~$\GF$ be a global field.  That is, $\GF$ is a finite extension of
the rational numbers or the field $\F_p(t)$, where~$\F_p$ is the
field with~$p$ elements for some prime~$p$.  The adele
ring~$\Adel_\GF$ over~$\GF$ is a certain locally compact ring that
contains~$\GF$ as a discrete cocompact subfield.  Its multiplicative
group~$\Adel\inv_\GF$ is a locally compact group as well for an
appropriate topology.  It contains~$\GF\inv$ as a discrete subgroup.
The quotient group $\ICL_\GF\defeq\Adel\inv_\GF/\GF\inv$ is the idele
class group of~$\GF$.

Let $\Sch(\Adel_\GF)$ be the space of Bruhat-Schwartz functions
on~$\Adel_\GF$ as defined by François Bruhat
in~\cite{Bruhat:Distributions}.  We equip it with a canonical
bornology, so that $\Sch(\Adel_\GF)$ is a nuclear, complete, convex
bornological vector space.  Both the pointwise product and the
convolution define algebra structures on $\Sch(\Adel_\GF)$.  There
exists an isomorphism between~$\Adel_\GF$ and its Pontrjagin
dual~$\widehat{\Adel_\GF}$.  Hence the Fourier transform can be viewed
as an operator on $\Sch(\Adel_\GF)$.  Since it intertwines the
convolution and the pointwise product, the two multiplications on
$\Sch(\Adel_\GF)$ give rise to isomorphic algebras.

The group~$\Adel\inv_\GF$ acts on~$\Adel_\GF$ by multiplication.
Thus~$\Adel\inv_\GF$ and its discrete subgroup~$\GF\inv$ act on
$\Sch(\Adel_\GF)$.  We let $\GF\inv\cross\Sch(\Adel_\GF)$ be the
algebraic crossed product.  As a bornological vector space, this is
$\bigoplus_{g\in\GF\inv} \Sch(\Adel_\GF)$ with the usual convolution
product.  We shall compute the Hochschild homology, the cyclic
homology and the periodic cyclic homology of this crossed product.
In addition, we also partially compute the topological
$\K$\nbd{}theory of the enveloping $C^*$\nbd{}algebra $\GF\inv\cross
C_0(\Adel_\GF)$.  We use~$\cross$ for both the algebraic and the
$C^*$\nbd{}algebra crossed products.  It should be clear from the
context which crossed product is meant.

An important additional structure on these homology spaces is the
action of the idele class group~$\ICL_\GF$.  Since~$\Adel\inv_\GF$ is
commutative, its action on $\Sch(\Adel_\GF)$ extends to an action on
$\GF\inv\cross\Sch(\Adel_\GF)$ by automorphisms.  Clearly, the
subgroup~$\GF\inv$ acts by inner automorphisms on the crossed product.
Therefore, it acts trivially on the various homology spaces of
$\GF\inv\cross\Sch(\Adel_\GF)$.  Thus the action of~$\Adel\inv_\GF$
descends to an action of~$\ICL_\GF$.  We also describe the
representation of~$\ICL_\GF$ on the homology spaces.  We state our
results and briefly summarize the way we obtain them in
Section~\ref{sec:results}.

Our computation of Hochschild and cyclic homology is a byproduct of
the spectral interpretation for poles and zeros of $L$\nbd{}functions
constructed by the author in~\cite{Meyer:Primes_Rep}.  This spectral
interpretation is inspired by a similar construction by Alain Connes
in~\cite{Connes:Trace_Formula}.  A basic ingredient is the coinvariant
space $\Sch(\Adel_\GF)/\GF\inv$.  However, to see zeros of
$L$\nbd{}functions, we have to compare this space to another natural
space $A(\ICL_\GF)$ of functions on~$\ICL_\GF$.  If~$\GF$ is a global
function field, then $A(\ICL_\GF)=\CCINF(\ICL_\GF)$ is the space of
locally constant functions with compact support.  If~$\GF$ is an
algebraic number field, then $A(\ICL_\GF)$ is the space of smooth
functions on~$\ICL_\GF$ with superexponentially decaying derivatives
for $\ln{}\abs{x}\to\pm\infty$.  Roughly speaking, we compare the
noncommutative quotient space $\Adel_\GF/\GF\inv$ to
$\Adel\inv_\GF/\GF\inv$.  Since the action of~$\GF\inv$
on~$\Adel\inv_\GF$ is free and proper, there is nothing noncommutative
about the latter quotient space.

The results of~\cite{Meyer:Primes_Rep} indicate that the spaces
$\Adel_\GF/\GF\inv$ and $\Adel\inv_\GF/\GF\inv$ are closely related.
This also comes out in our homology and $\K$\nbd{}theory computations.
The Hochschild, cyclic and periodic cyclic homology spaces of crossed
products decompose naturally as a direct sum of a homogeneous and an
inhomogeneous part.  We recall this decomposition later.  The
inhomogeneous part vanishes for $\Adel\inv_\GF/\GF\inv$.  It does not
vanish for $\Adel_\GF/\GF\inv$, but it is not very interesting either:
the inhomogeneous part can be detected by evaluation at the fixed
point $0\in\Adel_\GF$.  The homogeneous parts of the various cyclic
homology spaces for $\GF\inv\cross\Sch(\Adel_\GF)$ are very similar to
the corresponding spaces for $A(\ICL_\GF)$.  For the $\K$\nbd{}theory
computation, the relationship is not that strong but still
perceptible.

Éric Leichtnam and Victor Nistor have already attempted to compute the
cyclic theory of $\Q\inv\cross\Sch(\Adel_\Q)$
in~\cite{Leichtnam-Nistor:Adelic}.  They compute the inhomogeneous
part correctly, but they make a fatal mistake for the homogeneous
part, which leads to a wrong final result.  In particular, their
results for the Hochschild and cyclic homology of
$\Q\inv\cross\Sch(\Adel_\Q)$ have not much in common with the
corresponding spaces for $A(\ICL_\Q)$.  The general theory of how to
compute Hochschild and cyclic homology for group algebras and crossed
products can be found in~\cites{Burghelea, Nistor:Crossed,
Nistor:Crossed_algebraic, Getzler-Jones, Block-Getzler-Jones}.
However, we only use some elementary facts from this theory.  The
result of our Hochschild homology is so simple that we can pass to the
cyclic theories immediately.

Our $\K$\nbd{}theory computation uses the Baum-Connes conjecture for
the Abelian group~$\GF\inv$ and a spectral sequence for the
topological side of the Baum-Connes assembly map for torsion free
discrete groups.  Again, we use special properties of our situation to
simplify the computation.

\section{Summary and statement of the results}
\label{sec:results}

For any discrete group crossed product $G\cross A$, the Hochschild,
cyclic and periodic cyclic homology spaces split naturally over
conjugacy classes of~$G$.  The contribution of the trivial conjugacy
class of~$G$ is called the homogeneous part $(\cdots)_\homogen$, the
contribution of the other conjugacy classes is called the
inhomogeneous part $(\cdots)_\inhomogen$.

The homogeneous part of the Hochschild homology of
$\GF\inv\cross\Sch(\Adel_\GF)$ is easy to compute using results
of~\cite{Meyer:Primes_Rep}.  To describe it explicitly, we split
$\Adel_\GF=\Adel_\fin\times\Adel_\infty$, where $\Adel_\fin$
and~$\Adel_\infty$ are the finite and infinite adele rings,
respectively.  Thus~$\Adel_\fin$ is totally disconnected and
$\Adel_\infty\cong\R^n$ for some $n\in\N$.  Let
$\Lambda^n(\Adel^*_\infty)$ be the finite dimensional vector space of
$n$\nbd{}forms over the real vector space~$\Adel_\infty$ and let
$\Sch\Omega^n(\Adel_\GF)$ be the space of Bruhat-Schwartz functions
on~$\Adel_\GF$ taking values in $\Lambda^n (\Adel^*_\infty)$.  The
space $\Sch\Omega^n(\Adel_\GF)$ is naturally isomorphic to
$\HH_n(\Sch(\Adel_\GF))$.  Naturality includes the assertion that this
isomorphism is $\Adel\inv_\GF$\nbd{}equivariant.  Extending slightly
the results of~\cite{Meyer:Primes_Rep}, we show that
$H_j(\GF\inv,\Sch\Omega^n(\Adel_\GF))=0$ for $j\ge1$ and describe the
coinvariant space as a space of functions on~$\ICL_\GF$.  It follows
easily that
$$
\HH_n(\GF\inv\cross\Sch(\Adel_\GF))_\homogen \cong
\Sch\Omega^n(\Adel_\GF)/\GF\inv.
$$

This isomorphism on Hochschild homology can be realized by an explicit
chain map from the Hochschild complex
$(\Omega(\GF\inv\cross\Sch(\Adel_\GF)),b)$ to
$\Sch\Omega^*(\Adel_\GF)/\GF\inv$ with vanishing boundary map.  Up to
a constant factor, this chain map also intertwines the de Rham
boundary map $d_\dR\colon \Sch\Omega^n(\Adel_\GF) \to
\Sch\Omega^{n+1}(\Adel_\GF)$ and the boundary~$B$ on
$\Omega(\GF\inv\cross\Sch(\Adel_\GF))$.  Hence it induces a
quasi-isomorphism of the cyclic bicomplexes.  Since the vertical
boundary map in $\Sch\Omega^*(\Adel_\GF)$ vanishes, this allows us to
compute the cyclic homology and the periodic cyclic homology.

We need more notation in order to describe the result.  Let
$\sigma\colon \Adel\inv_\GF\to\{\pm1\}$ be the character that
multiplies the signs of an idele at the real places and let
$\C(\sigma)$ be the corresponding $1$\nbd{}dimensional representation.
There is a canonical integration map ${\smallint}\colon
\Sch\Omega^n(\Adel_\GF)\to\Sch(\Adel_\fin)$ that only depends on a
choice of orientation on~$\Adel_\fin$.  Hence it is an equivariant map
to the representation $\Sch(\Adel_\fin)\otimes\C(\sigma)$ of
$\Adel\inv_\GF\cong\Adel\inv_\fin\times\Adel\inv_\infty$.  It is easy
to see that the kernel of~$\smallint$ is a contractible complex.
Hence $(\Sch\Omega^*(\Adel_\GF),d_\dR,{\smallint})$ is a resolution of
$\Sch(\Adel_\fin)\otimes\C(\sigma)$.  Our previous results show that
this resolution is acyclic for group homology.  Hence the complex
$(\Sch\Omega^*(\Adel_\GF)/\GF\inv,d_\dR)$ computes the group homology
$H_j(\GF\inv,\Sch(\Adel_\fin)\otimes\C(\sigma))$.  We obtain
$$
\HP_j(\GF\inv\cross \Sch(\Adel_\fin))_\homogen \cong
\bigoplus_{k\in\Z}
H_{j+2k}(\GF\inv,\Sch(\Adel_\fin)\otimes\C(\sigma))
$$
for $j\in\Ztwo$ and a similar result for cyclic homology.  It remains
to describe $H_j(\GF\inv,\Sch(\Adel_\fin)\otimes\C(\sigma))$.  The
higher group homology spaces vanish for $\GF=\Q$ and more generally
if~$\GF$ has at most one infinite place.  Let $\GF_\infty\subseteq\GF$
be the subring of algebraic integers and let $\GF\inv_\infty$ be the
multiplicative group of invertible algebraic
integers~$\GF\inv_\infty$.  It is known that~$\GF\inv_\infty$ is a
product of the finite cyclic group of roots of unity in~$\GF\inv$ and
a free Abelian group of rank $u-1$, where~$u$ is the number of
infinite places of~$\GF$.  The action of~$\GF\inv_\infty$ on
$\Sch(\Adel_\fin)\hot\C(\sigma)$ factors through a compact group and
hence is close to a trivial action.  Hence it is no surprise that the
group homology of $\Z^{u-1}$ comes up.  Indeed, there is an
isomorphism
$$
H_k(\GF\inv,\Sch(\Adel_\fin)\otimes\C(\sigma)) \cong
\CCINF(\Adel\inv_\fin/\cl{\GF\inv}) \hot \Lambda^k \C^{u-1},
$$
where $\cl{\GF\inv}$ denotes the closure of~$\GF\inv$ in
$\Adel\inv_\fin$.

The inhomogeneous part of the Hochschild homology can be described as
follows.  Evaluation at the fixed point $0\in\Adel_\GF$ of the action
of~$\Adel\inv_\GF$ induces an algebra homomorphism
$\GF\inv\cross\Sch(\Adel_\GF)\to\C[\GF\inv]$.  This homomorphism
induces an isomorphism on the inhomogeneous part of the Hochschild
homology.  Hence it also induces an isomorphism on the inhomogeneous
parts of cyclic and periodic cyclic homology.

It is straightforward to compute the various homology spaces for
$\C[\GF\inv]$, using that~$\GF\inv$ is isomorphic to a product of the
finite cyclic group~$\Roots$ of roots of unity in~$\GF\inv$ and a free
Abelian group of infinite rank.  Hence the Pontrjagin
dual~$\widehat{\GF\inv}$ is a product of the finite cyclic
group~$\widehat{\Roots}$ and an infinite dimensional torus.  The
Fourier transform identifes $\C[\GF\inv]$ with the algebra of Laurent
series on~$\widehat{\GF\inv}$.  Hence the $n$th Hochschild homology is
the space $\Omega^n(\widehat{\GF\inv})$ of differential $n$\nbd{}forms
on~$\widehat{\GF\inv}$ with Laurent series as coefficients.  The
periodic cyclic homology is the de Rham cohomology
$H_\dR^*(\widehat{\GF\inv})$ of~$\widehat{\GF\inv}$.  The latter is a
direct sum, indexed by~$\widehat{\Roots}$, of copies of the exterior
algebra in infinitely many generators.

Next we explain the $\K$\nbd{}theory computation.  First, we have to
compute the equivariant $\K$\nbd{}theory group
$\K^\Roots(C_0(\Adel_\GF))$, where~$\Roots$ is the group of roots of
unity in~$\GF\inv$ as above.  We obtain a natural splitting
$$
\K^\Roots(C_0(\Adel_\GF))\cong \K_\homogen \oplus \K_\inhomogen
$$
that corresponds to the homogeneous and inhomogeneous parts of the
homology computation.  The summand $\K_\inhomogen$ is again detected
by evaluation at $0\in\Adel_\GF$.  The summand $\K_\homogen$ is
isomorphic to the space of compactly supported locally constant
functions $f\colon \Adel_\fin\to\Z$ that satisfy $f(ax)=\sigma(a)f(x)$
for all $a\in\Roots$, where $\sigma(a)=\pm1$.  Thus we obtain an
integral lattice in the representation
$\Sch(\Adel_\fin)\otimes\C(\sigma)$ that occured in the computation of
the periodic cyclic homology.

We can write
$$
\GF\inv\cross C_0(\Adel_\GF) \cong
\GF\inv/\Roots \cross (\Roots\cross C_0(\Adel_\GF)).
$$
The group $\GF\inv/\Roots$ is torsion free and satisfies the
Baum-Connes conjecture with coefficients.  Hence $\K(\GF\inv\cross
C_0(\Adel_\GF))$ is the limit of a spectral sequence with
$$
E^2_{pq} =
H_p(\GF\inv/\Roots,\K^\Roots_q(C_0(\Adel_\GF))) =
H_p(\GF\inv/\Roots,\K_\homogen) \oplus
H_p(\GF\inv/\Roots,\K_\inhomogen).
$$
We check that the spectral sequence already collapses at~$E^2$.  Since
the group~$\ICL_\GF$ acts differently on $\K_\homogen$ and
$\K_\inhomogen$, there is no interaction between the corresponding
parts of the spectral sequence.  Hence we only have to consider the
homogeneous part of $E^2_{pq}$.  It describes the subquotients of a
certain filtration on the ``homogeneous part'' of the
$\K$\nbd{}theory.  If $\GF=\Q$ or more generally if~$\GF$ has just one
infinite place, then $E^2_{pq}=0$ for $p\ge1$, so that we obtain the
$\K$\nbd{}theory itself.  In general, we should also describe the
extensions by which we get the $\K$\nbd{}theory from the subquotients
$E^2_{pq}$.  We do not do that.

Our computation of $H_p(\GF\inv/\Roots,\K_\homogen)$ proceeds as
follows.  We first restrict the action to $\GF\inv_\infty/\Roots
\subseteq \GF\inv/\Roots$.  This yields a space of functions
from~$\Adel_\fin$ to an exterior algebra as for the periodic cyclic
homology.  However, these functions are allowed to take certain
rational values, with a rather complicated rule for the denominators
that are allowed at a point $x\in\Adel_\fin$.  Nevertheless, we can
show that
$$
H_j(\GF\inv/\GF\inv_\infty, H_k(\GF\inv_\infty/\Roots, \K_\homogen))
\cong 0
$$
for $j>0$, so that
$$
E^2_{pq} \cong
H_0(\GF\inv/\GF\inv_\infty, H_p(\GF\inv_\infty/\Roots, \K_\homogen))
$$
The proof uses the filtrations on $\Sch(\Adel_\fin,\Z)$ and
$\CCINF(\Adel\inv_\fin,\Z)$ by the subspaces of functions that are
ramified at at most~$n$ places.  The subquotients of these filtrations
on $\Sch(\Adel_\fin,\Z)$ and $\CCINF(\Adel\inv_\fin,\Z)$ turn out to
be isomorphic and can be described explicitly.

\section{Preparations}
\label{sec:preparations}

We begin by recalling our notation regarding adeles and ideles and
Bruhat-Schwartz spaces.  We are going to use some general results
about homological algebra for certain smooth convolution algebras on
locally compact groups that are proved already
in~\cite{Meyer:Primes_Rep}.  We identify
$\GF\inv\cross\Sch(\Adel_\GF)$ with a smooth convolution algebra on
the crossed product group $\GF\inv\cross\Adel_\GF$ and show that this
convolution algebra is \emph{isocohomological}.  This allows us to
apply the usual recipe for computing the Hochschild homology of group
algebras.

\subsection{Adeles and Bruhat-Schwartz functions}
\label{sec:adeles_Sch}

We begin by recalling some basic notation and results concerning local
and global fields and the space $\Sch(\Adel_\GF)$.  More details can
be found in~\cites{Bruhat:Distributions, Meyer:Primes_Rep,
Weil:Basic}.

\emph{Local fields} can be described abstractly as non-discrete,
locally compact fields.  More concretely, the local fields can be
defined as finite extensions of the fields $\R$, $\Q_p$ and
$\F_p((t))$ for prime numbers~$p$.  Here $\F_p((t))$ is obtained by
adjoining~$t^{-1}$ to the ring $\F_p[[t]]$ of formal power series
over~$\F_p$ and~$\Q_p$ denotes as usual the $p$\nbd{}adic integers.
Of course, the only finite extensions of~$\R$ are $\R$ and~$\C$, these
are the Archimedean or \emph{infinite} local fields.  The other local
fields are called \emph{finite} or non-Archimedean.  The only finite
extensions of~$\F_p((t))$ are $\F_q((t))$, where~$q$ is a power
of~$p$.  We cannot list the finite extensions of~$\Q_p$ so easily.

Each local field is equipped with a canonical \emph{norm
homomorphism}, $x\mapsto\abs{x}$.  This is the usual $p$\nbd{}adic
norm~$\abs{x}_p$ for~$\Q_p$; the usual absolute value for~$\R$; and
the \emph{square} of the usual absolute value for~$\C$.

A \emph{global field} is, by definition, a finite extension of~$\Q$
or~$\F_p(t)$ for some prime~$p$.  Global fields of the first kind are
called \emph{algebraic number fields}, those of the second kind are
called \emph{global function fields}.  A \emph{place} of a global
field~$\GF$ is a dense embedding of~$\GF$ into a local field.  For
instance, the places of~$\Q$ are the embeddings $\Q\to\Q_p$ for the
prime numbers~$p$ and the single infinite place $\Q\to\R$.  A place is
called finite, infinite, real or complex depending on the local
field~$\GF_v$.

We denote the set of places of~$\GF$ by $\Places(\GF)$.  Let
$S\subseteq\Places(\GF)$ be a subset that contains all infinite places
of~$\GF$ and let $\Cont{S}\defeq \Places(\GF)\setminus S$.  We let
\begin{align*}
  \GF_S &\defeq
  \{a\in\GF\mid \text{$\abs{a}_v\le1$ for all $v\in\Cont{S}$}\},
  \\
  \GF\inv_S &\defeq
  \{a\in\GF\mid \text{$\abs{a}_v=1$ for all $v\in\Cont{S}$}\}.
\end{align*}
Then $\GF_S\subseteq\GF$ is a subring and~$\GF\inv_S$ is its
multiplicative group of units.  That is, $a\in\GF\inv_S$ if and only
if both $a$ and~$a^{-1}$ belong to~$\GF_S$.  If~$\GF$ is an algebraic
number field and~$S$ is the set of all infinite places, then
$\GF_S\subseteq\GF$ is equal to the ring of \emph{algebraic integers}
in~$\GF$.  We also denote this subring by~$\GF_\infty$.  If~$\GF$ is a
global function field, then there are no infinite places, so that
$S=\emptyset$ is admissible.  We let $\GF_\infty\defeq\GF_\emptyset$
in this case.

The \emph{adele ring}~$\Adel_\GF$ of~$\GF$ is the \emph{restricted}
direct product of the local fields~$\GF_v$, where~$v$ runs through the
places of~$\GF$.  This is a locally compact ring, that is, its
topology is locally compact and addition and multiplication are
continuous.  The embeddings $\GF\to\GF_v$ combine to a \emph{diagonal
embedding} $\GF\to\Adel_\GF$.  Its range is a discrete and cocompact
subfield of~$\Adel_\GF$.  The \emph{idele group}~$\Adel\inv_\GF$ is
the multiplicative group of~$\Adel\inv_\GF$.  It can also be described
as the restricted direct product of the multiplicative
groups~$\GF\inv_v$ for $v\in\Places(\GF)$ and is a locally compact
Abelian group as well.  The \emph{idele class group} is the quotient
group $\ICL_\GF\defeq \Adel\inv_\GF/\GF\inv$.  Since
$\GF\inv\subseteq\Adel\inv_\GF$ is discrete, this is again a locally
compact group.

The norm homomorpisms on the groups~$\GF\inv_v$ combine to a norm
homomorphism
$$
\abs{\blank}\colon
\Adel\inv_\GF\to\R\inv_+,
\qquad x\mapsto\abs{x},
$$
which is trivial on~$\GF\inv$ and hence descends to the quotient
group~$\ICL_\GF$.

Given a set~$S$ of places of~$\GF$, we can also define the
\emph{$S$\nbd{}adele ring}~$\Adel_S$ and the \emph{$S$\nbd{}idele
group}~$\Adel\inv_S$ by taking the restricted direct products of
$\GF_v$ and~$\GF\inv_v$ only over $v\in S$.  We are particularly
interested in the cases where~$S$ consists of all finite or of all
infinite places of~$\GF$.  These give rise to the rings $\Adel_\fin$
and~$\Adel_\infty$ of finite and infinite adeles and the groups
$\Adel\inv_\fin$ and~$\Adel\inv_\infty$ of finite and infinite ideles.
Let $n_\GF=[\GF:\Q]$ for algebraic number fields and $n_\GF=0$ for
global function fields.  Then $\Adel_\infty\cong\R^{n_\GF}$ as a
locally compact Abelian group.

The space $\Sch(\Adel_\GF)$ is defined by François Bruhat
in~\cite{Bruhat:Distributions}.  We view it as a bornological vector
space as in~\cite{Meyer:Primes_Rep}.  The product decomposition
$\Adel_\GF=\Adel_\fin\times\Adel_\infty$ gives rise to a tensor
product decomposition
$$
\Sch(\Adel_\GF)\cong\Sch(\Adel_\fin)\hot\Sch(\Adel_\infty).
$$
Since~$\Adel_\fin$ is totally disconnected and a direct union of
compact subgroups, $\Sch(\Adel_\fin)$ is equal to the space of locally
constant, compactly supported functions on~$\Adel_\fin$.  This is a
direct union of finite dimensional subspaces and hence carries the
fine bornology.  Since $\Adel_\infty\cong\R^n$ for some $n\in\N$,
$\Sch(\Adel_\infty)$ is the usual Schwartz space on~$\R^n$.  A subset
is bounded if and only if it is von Neumann bounded in the usual
Fréchet topology on $\Sch(\R^n)$.  Since $\Sch(\Adel_\fin)$ carries
the fine bornology, we have $\Sch(\Adel_\GF)\cong
\Sch(\Adel_\fin)\otimes\Sch(\Adel_\infty)$, that is, the algebraic
tensor product is already complete.

Both the pointwise multiplication~$\cdot$ and the convolution~$*$
(with respect to the additive structure) turn $\Sch(\Adel_\GF)$ into
bornological algebras.  The Fourier transform gives rise to an
isomorphism between these two algebras, as follows.  Since
$\GF\subseteq\Adel_\GF$ is discrete, there is a non-trivial character
$\psi\colon\Adel_\GF\to S^1$ with $\psi(a)=1$ for all $a\in\GF$.  Then
the bilinear pairing
$$
\Adel_\GF\times\Adel_\GF\to S^1,
\qquad
(x,y)\mapsto \psi(x\cdot y)
$$
induces an isomorphism of locally compact Abelian groups
$\Adel_\GF\cong\widehat{\Adel_\GF}$.  We use this isomorphism to
define the Fourier transform
\begin{equation}  \label{eq:def_Fourier}
  \Fourier\colon \Sch(\Adel_\GF)\to \Sch(\Adel_\GF),
  \qquad
  \Fourier f(\xi)\defeq
  \int_{\Adel_\GF} f(x) \psi(x\xi)\,dx.
\end{equation}
We normalize the Haar measure~$dx$ on~$\Adel_\GF$ so that
$\Fourier^* f(x)=\Fourier f(-x)$ is inverse to~$\Fourier$.
The Fourier transform is an isomorphism of bornological algebras
$$
(\Sch(\Adel_\GF),*) \cong (\Sch(\Adel_\GF),\cdot).
$$

The group~$\Adel\inv_\GF$ acts on~$\Adel_\GF$ by multiplication.  Let
$$
\lambda_g f(x)\defeq f(g^{-1}x),
\qquad
\rho_g f(x)\defeq \abs{g}\cdot f(xg),
$$
for $g\in\Adel\inv_\GF$, $f\in\Sch(\Adel_\GF)$, $x\in\Adel_\GF$.
Then~$\lambda$ and~$\rho$ are smooth representations
of~$\Adel\inv_\GF$ on $\Sch(\Adel_\GF)$.  The maps $\lambda_g$
and~$\rho_g$ are algebra automorphisms with respect to pointwise
multiplication and convolution, respectively.  Since $\Fourier\circ
\lambda_g = \rho_g \circ \Fourier$, we have an isomorphism of
dynamical systems
$$
(\Sch(\Adel_\GF),\cdot,\lambda)\cong
(\Sch(\Adel_\GF),*,\rho),
$$
so that we need not distinguish between them.

We define $\GF\inv\cross \Sch(\Adel_\GF)$ to be
$\bigoplus_{g\in\GF\inv} \Sch(\Adel_\GF)$ equipped with the direct sum
bornology and the usual multiplication of a crossed product.

\subsection{Isocohomological smooth convolution algebras}
\label{sec:isocohomological}

There is a semidirect product \emph{group} $G\defeq
\GF\inv\cross\Adel_\GF$ with multiplication
$$
G\times G\to G,
\qquad
(g_1,x_1)\cdot (g_2,x_2)=(g_1g_2,g_2^{-1}x_1+x_2).
$$
The product in the crossed product
$$
A(G)\defeq \GF\inv\cross(\Sch(\Adel_\GF),*)
$$
is exactly the convolution for this group structure on~$G$.  It is
clear that $A(G)$ contains $\CCINF(G)$ as a dense subalgebra.
The left and right regular representations of~$G$ on $A(G)$ are
smooth because the regular representation of~$\Adel_\GF$ on
$\Sch(\Adel_\GF)$ is smooth.  Thus
$A(G)$ is a smooth convolution algebra on~$G$ in the sense
of \cite{Meyer:Primes_Rep}*{Section 4.3}.

We define the category $\Mod_{A(G)}$ of essential (left) modules over
$A(G)$ as in \cite{Meyer:Primes_Rep}.  We use the notation
$\Mod_{A(G)}$ and ${}_{A(G)}\Mod$ for the categories of essential
right and left $A(G)$\nbd{}modules if we need both kinds of modules.
We usually avoid the ugly notation ${}_{A(G)}\Mod$.  Recall that
$\Mod_{A(G)}$ is isomorphic to a full subcategory of the category of
smooth representations of~$G$.  Let $A(G)^\op$ be the opposite algebra
of~$A(G)$.  Then the category of essential $A(G)$\nbd{}bimodules is
isomorphic to the category of essential modules over $A(G)\hot
A(G)^\op$.  The latter is a smooth convolution algebra on $G\times
G^\op\cong G\times G$.  The quickest way to see this is to observe
that both categories of essential modules are isomorphic to the same
category of smooth representations.

We require extensions in $\Mod_{A(G)}$ to split as extensions of
bornological vector spaces.  This turns $\Mod_{A(G)}$ into an exact
category, so that we can form its derived category $\Der_{A(G)}$,
whose objects are chain complexes over $\Mod_{A(G)}$.  Since modules
over the zero algebra are nothing but (complete convex) bornological
vector spaces, the category $\Der_0$ is the homotopy category of chain
complexes of bornological vector spaces.  In the following, all
bornological vector spaces are tacitly assumed complete and convex.
We denote the total left derived functor of the balanced tensor
product functor $\hot_{A(G)}$ by $\Lhot_{A(G)}$.  This is a functor
$$
\Lhot_{A(G)}\colon
\Der_{A(G)} \times {}_{A(G)}\Der \to \Der_0.
$$
We obtain the more classical satellite functors by embedding
$\Mod_{A(G)}\to\Der_{A(G)}$ in the usual way and taking the homology
of the chain complex $V\Lhot_{A(G)} W$.  It is convenient for
statements to replace the usual Hochschild homology by a chain complex
as well.  We let $\CHH(A(G))$ be the image of the chain complex
$(\Omega^n A(G),b)$ in $\Der_0$.  Thus the Hochschild homology of
$A(G)$ is the homology of $\CHH(A(G))$.

It is observed in~\cite{Meyer:Primes_Rep} that the algebra $A(G)$
contains an approximate identity and that $A(G)$ is projective both as
a left and right $A(G)$\nbd{}module.  This is true for any smooth
convolution algebra and follows easily from the corresponding results
for $\CCINF(G)$ proved in~\cite{Meyer:Smoothrep}.  It follows that any
essential module has a (linearly split) resolution by essential
modules that are projective among not necessarily essential
$A(G)$\nbd{}modules.  Thus it makes no difference whether we compute
derived functor in the category of all modules or the category of
essential modules.  It also follows that $A(G)$ is $H$\nbd{}unital.
Hence we can rewrite Hochschild homology as a balanced tensor product
functor
$$
\CHH(A(G)) \cong A(G) \Lhot_{A(G)\hot A(G)^\op} A(G).
$$

Let $\C(1)$ be the trivial representation of~$G$, viewed as an
essential right module over $\CCINF(G)$, and let $\CCINF(G)_\Ad$ and
$A(G)_\Ad$ denote the representation of~$G$ on $\CCINF(G)$ and $A(G)$
by conjugation, $g\cdot f(x)\defeq f(g^{-1}xg)$.  If~$W$ is a smooth
representation of~$G$, we let
$$
\CHo(G,W) \defeq \C(1)\Lhot_{\CCINF(G)} W \in \Der_0.
$$
By definition, the homology of this chain complex is the \emph{group
homology} of~$G$.  We would like to have an isomorphism
$\CHH(A(G))\cong \CHo(G,A(G)_\Ad)$.  However, this is not true
for arbitrary smooth convolution algebras.  It fails, for instance,
for $\ell^1(\Z)$, which is a smooth convolution algebra on~$\Z$.  We
only get such an isomorphism if $A(G)$ is \emph{isocohomological} in
the sense of~\cite{Meyer:Primes_Rep}.

A smooth convolution algebra $A(G)$ on a locally compact group is
called \emph{isocohomological} if the fully faithful embedding
$\Mod_{A(G)}\to\Mod_{\CCINF(G)}$ gives rise to a fully faithful
functor $\Der_{A(G)}\to\Der_{\CCINF(G)}$.  An equivalent condition is
that $V\Lhot_{A(G)} W\cong V\Lhot_{\CCINF(G)} W$ for all
$V\in\Der_{A(G)}$, $W\in{}_{A(G)}\Der$.  Even more, it suffices to
check
$$
A(G) \Lhot_{\CCINF(G)} A(G) \cong
A(G) \Lhot_{A(G)} A(G) \cong A(G)
$$
because this implies the corresponding assertion for chain complexes
of free essential $A(G)$\nbd{}modules.

\begin{lemma}  \label{lem:our_crossed_product_isocohomological}
  The crossed product $\GF\inv\cross \Sch(\Adel_\GF)$ is an
  isocohomological smooth convolution algebra on
  $\GF\inv\cross\Adel_\GF$.
\end{lemma}

\begin{proof}
  Let $H\defeq\Adel_\GF$ and $G=\GF\inv\cross\Adel_\GF$ and let
  $A(G)\defeq \GF\inv\cross \Sch(\Adel_\GF)$.  Then~$H$ is an open
  normal subgroup of~$G$ and there is a natural isomorphism
  $$
  A(G) \cong \cInd_H^G \Sch(H)
  $$
  as left or right modules over $\CCINF(G)$.  We compute
  \begin{multline*}
    A(G) \Lhot_{\CCINF(G)} A(G)
    \cong
    \cInd_H^G \Sch(H) \Lhot_{\CCINF(G)} A(G)
    \cong
    \Sch(H) \Lhot_{\CCINF(H)} \Res_G^H A(G)
    \\ \cong
    \Sch(H) \Lhot_{\Sch(H)} \Res_G^H A(G)
    \cong
    \Res_G^H A(G),
  \end{multline*}
  using that compact induction and restriction are adjoint
  (\cite{Meyer:Smoothrep}) and that $\Sch(\Adel_\GF)$ is
  isocohomological (\cite{Meyer:Primes_Rep}).
\end{proof}

\begin{proposition}  \label{pro:isocohomological_HH}
  Let $A(G)$ be an isocohomological smooth convolution algebra on a
  locally compact group~$G$.  Then
  $$
  \CHH(A(G)) \cong \CHo(G,A(G)_\Ad).
  $$
\end{proposition}

\begin{proof}
  The condition $A(G)\Lhot_{\CCINF(G)} A(G)\cong A(G)$ is symmetric
  under passage to opposite algebras and hereditary for tensor
  products.  Hence $A(G)\hot A(G)^\op$ is an isocohomological smooth
  convolution algebra on $G\times G$.  Therefore,
  $$
  \CHH(A(G)) \cong
  A(G) \Lhot_{A(G)\hot A(G)^\op} A(G) \cong
  A(G) \Lhot_{\CCINF(G\times G)} A(G).
  $$
  It is shown in~\cite{Meyer:Smoothrep} that the bivariant functor
  $\Lhot_{\CCINF(G\times G)}$ can be reduced to group homology.  In
  our case, this yields the following.  Let $G\times G$ act on
  $A(G)\hot A(G)$ by 
  $$
  (g_1,g_2) \bullet f(h_1,h_2) = f(g_1^{-1}h_1g_2,g_2^{-1}h_2g_1).
  $$
  Then
  \begin{multline*}
    A(G) \Lhot_{\CCINF(G\times G)} A(G) \cong
    \C(1) \Lhot_{\CCINF(G\times G)} A(G)\hot A(G) \\ \cong
    \C(1) \Lhot_{\CCINF(G_1)} (\C(1) \Lhot_{\CCINF(G_2)}
    A(G)\hot A(G)).
  \end{multline*}
  Here $G_1$ and~$G_2$ denote the first and second factor~$G$ in
  $G\times G$.  Reversing the above argument, we can rewrite
  $$
  \C(1) \Lhot_{\CCINF(G_2)} A(G)\hot A(G) \cong
  A(G) \Lhot_{\CCINF(G_2)} A(G) \cong A(G),
  $$
  where we used once again that~$A(G)$ is isocohomological.  The
  resulting action of $G=G_1$ on $A(G)$ is by conjugation, so that we
  obtain $\C(1)\Lhot_{\CCINF(G)} A(G)_\Ad$ as desired.
\end{proof}

This proposition applies to $G\defeq\GF\inv\cross\Adel_\GF$ and
$A(G)\defeq\GF\inv\cross \Sch(\Adel_\GF)$ by
Lemma~\ref{lem:our_crossed_product_isocohomological}.  Hence we obtain
$$
\CHH(\GF\inv\cross\Sch(\Adel_\GF)) =
\CHH(A(G)) \cong
\CHo(G, A(G)_\Ad).
$$

\subsection{The coaction on the homology spaces}
\label{sec:coaction_GFinv}

Next we explain how the coaction of~$\GF\inv$ on the various cyclic
homology groups comes about.  Since~$\GF\inv$ is discrete, a coaction
on a bornological vector space~$V$ is nothing but a direct sum
decomposition $V=\bigoplus_{g\in\GF\inv} V_g$.  For coactions on
algebras, we require $A_gA_h\subseteq A_{gh}$ for all $g,h\in G$.
If~$A$ is any algebra on which~$\GF\inv$ acts by automorphisms, then
$\GF\inv\cross A\defeq \sum_{g\in G} A_g$ with $A_g\defeq
A\delta_g\subseteq \GF\inv\cross A$ carries such a coaction.  We
describe it by the partially defined function $\deg(a)\defeq g$ for
$a\in A_g$.  The complex $\Omega (\GF\inv\cross A)$ of non-commutative
differential forms inherits a ``total'' grading
$$
\deg(a_0\,da_1\ldots\,da_n)\defeq \deg(a_0)\cdots \deg(a_n).
$$
The boundary maps $b$ and~$d$ preserves this grading because~$\GF\inv$
is commutative.  Hence we obtain a direct sum decomposition
$$
\CHH(\GF\inv\cross A) \cong
\bigoplus_{g\in\GF\inv} \CHH(\GF\inv\cross A)_g
$$
and a similar decomposition of Hochschild homology.  Since~$B$ also
preserves the coaction, we also have such decompositions for cyclic
and periodic cyclic homology.  It turns out that the
$1$\nbd{}homogeneous parts play a special role.  Hence we only
distinguish the \emph{homogeneous part} $(\cdots)_1$ and the
\emph{inhomogeneous part} $\bigoplus_{g\neq1} (\cdots)_g$

\section{Computation of the Hochschild homology}
\label{sec:HH}

As in Section~\ref{sec:isocohomological}, we let $G\defeq
\GF\inv\cross\Adel_\GF$ and view
$$
A(G)\defeq \GF\inv\cross(\Sch(\Adel_\GF),*)
$$
as a smooth convolution algebra on~$G$.  We have seen there that the
complex $\CHH(A(G))$ that computes the Hochschild homology of $A(G)$
is homotopy equivalent to the complex $\CHo(G,A(G)_\Ad)$ that computes
the group homology of~$G$ with coefficients in $A(G)$ with~$G$ acting
by conjugation.  Since~$G$ is defined as a semidirect product, we can
compute this group homology in two stages as
$$
\CHo(G,A(G)_\Ad) \cong
\CHo(\GF\inv,\CHo(\Adel_\GF,A(G)_\Ad)),
$$
where we view $\CHo(\Adel_\GF,A(G)_\Ad)$ as an object of
$\Der_{\CCINF(\GF\inv)}$ in the canonical way.  If we worked with
satellite functors, the above isomorphism would become the spectral
sequence computing group homology for a group extension.

The subsets $\{g\}\times\Adel_\GF\subseteq G$ are invariant under
conjugation because~$\GF\inv$ is commutative.  Hence
$$
A(G)_\Ad \cong \bigoplus_{g\in\GF\inv} \Sch(\Adel_\GF)_g,
$$
where $\Sch(\Adel_\GF)_g$ is $\Sch(\Adel_\GF)$ as a bornological
vector space.  The action of~$G$ on $\Sch(\Adel_\GF)$ comes from the
action on the subset $\{g\}\times\Adel_\GF$ by conjugation.  Since
$$
(h,y)^{-1}(g,x)(h,y) =
(h^{-1},-hy)(g,x)(h,y) =
(g,(1-g^{-1})y+h^{-1}x),
$$
the representation of~$G$ on $\Sch(\Adel_\GF)_g$ is given by
$$
(h,y) \bullet_g f(x) = f(h^{-1}x+ (1-g^{-1})y)
$$
for all $g,h\in\GF\inv$, $x,y\in\Adel_\GF$.

We have to distinguish the cases $g\neq1$ and $g=1$.  We treat the
case $g\neq1$ first.  Then $1-g^{-1}$ is invertible, so that the map
$\Phi f(x)=f\bigl((1-g^{-1})x\bigr)$ is a bornological isomorphism on
$\Sch(\Adel_\GF)_g$.  It intertwines the representation~$\bullet_g$ on
$\Sch(\Adel_\GF)_g$ described above and the standard representation
$$
(h,y)\bullet f(x) = f(h^{-1}x+y).
$$
Its restriction to $\Adel_\GF\subseteq G$ is the regular
representation.  Since $\C(1)$ is an essential module over
$\Sch(\Adel_\GF)$ and $\Sch(\Adel_\GF)$ is isocohomological
by~\cite{Meyer:Primes_Rep}, we obtain
$$
\CHo(\Adel_\GF,\Sch(\Adel_\GF)_g) =
\C(1)\Lhot_{\CCINF(\Adel_\GF)} \Sch(\Adel_\GF)_g \cong
\C(1)\Lhot_{\Sch(\Adel_\GF)} \Sch(\Adel_\GF) \cong
\C(1).
$$
This homotopy equivalence is implemented by the $G$\nbd{}invariant
linear map
$$
\Sch(\Adel_\GF)_g\to\C(1),
\qquad
f\mapsto\int_{\Adel_\GF} f(x)\,dx
$$
and hence is an isomorphism in the derived category
$\Der_{\CCINF(\GF\inv)}$.  Therefore,
$$
\CHH(\GF\inv\cross\Sch(\Adel_\GF))_g \cong
\CHo(\GF\inv,\CHo(\Adel_\GF,\Sch(\Adel_\GF)_g)) \cong
\CHo(\GF\inv,\C(1)).
$$
The homology of $\CHo(\GF\inv,\C(1))$ is the group homology of the
discrete group~$\GF\inv$ in the usual sense.  We will describe it more
concretely below.

Now we consider the homogeneous part, that is, we let $g=1$.  The
subgroup~$\Adel_\GF$ acts trivially on $\Sch(\Adel_\GF)_1$.  Therefore,
$$
\CHo(\Adel_\GF,\Sch(\Adel_\GF)_1) \cong
\CHo(\Adel_\GF,\C(1))\hot \Sch(\Adel_\GF),
$$
where $\CHo(\Adel_\GF,\C(1))$ is a chain complex of bornological
vector spaces that computes the group homology of~$\Adel_\GF$ with
trivial coefficients.  The direct product decomposition
$\Adel_\GF\cong\Adel_\fin\times\Adel_\infty$ gives rise to a tensor
product decomposition for $\CCINF(\Adel_\GF)$ and hence to an
isomorphism
$$
\CHo(\Adel_\GF,\C(1))\cong
\CHo(\Adel_\fin,\C(1))\hot\CHo(\Adel_\infty,\C(1)).
$$

We have already observed that $\CCINF(\Adel_\fin)\cong
\Sch(\Adel_\fin)$.  Hence the Fourier transform identifies
$(\CCINF(\Adel_\fin),*)$ and $(\CCINF(\Adel_\fin),\cdot)$.  The
space~$\Adel_\fin$ is totally disconnected.  Therefore, any essential
module over $\CCINF(\Adel_\fin)$ is flat, that is,
$V\Lhot_{\CCINF(\Adel_\fin)}W=V\hot_{\CCINF(\Adel_\fin)}W$ for all
$V,W$.  In particular, $\C(1) \Lhot_{\CCINF(\Adel_\fin)} \C(1) \cong
\C(1)$.  Since $\Adel_\infty\cong\R^{n_\GF}$, we have
$$
\CHo(\Adel_\infty,\C(1))\cong
\CHo(\R,\C(1))^{\hot n_\GF} \cong
(\Lambda \R)^{\hot n_\GF} \cong
\Lambda(\R^{n_\GF}).
$$
Here~$\Lambda(V)$ denotes the exterior algebra over a real vector
space~$V$.  We have to describe the action of~$\GF\inv$ on this space.
Therefore, we write down the above computation carefully and
naturally.

The space~$\Adel_\infty$ carries a canonical real vector space
structure.  The action of~$\Adel\inv_\infty$ on~$\Adel_\infty$ is
linear.  Let $\Lambda \Adel_\infty$ be the exterior algebra over the
real vector space~$\Adel_\infty$.  This carries an induced
representation of~$\Adel\inv_\infty$.  We let~$\Adel\inv_\fin$
act trivially on $\Adel_\infty$ and $\Lambda\Adel_\infty$ and denote
the resulting representation of~$\Adel\inv_\GF$ on
$\Lambda\Adel_\infty$ by~$\pi$.  Consider the graded vector space
$$
C(\Adel_\infty) \defeq \CCINF(\Adel_\infty)\hot \Lambda \Adel_\infty
$$
with the boundary map
$$
b(f\otimes \xi_1\wedge\dots\wedge \xi_j) \defeq
\sum_{k=1}^j (-1)^{k-1} \frac{\partial f}{\partial \xi_k}\otimes
\xi_1\wedge\dots \wedge\widehat{\xi_k}
\wedge\dots \wedge\xi_j
$$
of degree~$-1$ and the augmentation map
$$
{\smallint}\colon C_0(\Adel_\infty)\cong \CCINF(\Adel_\GF)\to\C,
\qquad
f\mapsto \int_{\Adel_\GF} f(x)\,dx.
$$
We have ${\smallint}\circ b=0$ and $b^2=0$ because the differential
operators $\partial/\partial\xi$ for $\xi\in\Adel_\infty$ commute.
Thus $(C(\Adel_\infty),b,{\smallint})$ is a bornological chain complex
over~$\C$.  Writing $\Adel_\infty\cong\R^n$, we get
$(C(\Adel_\infty),b) \cong(C(\R),b)^{\hot n}$.  It is easy to check
that $(C(\R),b,{\smallint})$ is a linearly split resolution of~$\C$.
Hence so is $(C(\Adel_\infty),b,{\smallint})$.

We let $\Adel\inv_\GF$ and~$\Adel_\GF$ act on $C(\Adel_\infty)$ by
$\rho\hot\pi$ and by translation on the first tensor factor.  Thus
$\Adel\inv_\fin\cross \Adel_\fin\subseteq\Adel_\GF$ acts trivially.
These representation piece together to a smooth representation of
$\Adel\inv_\GF\cross\Adel_\GF$.  The boundary map~$b$ is equivariant
for this representation and~${\smallint}$ is an invariant linear
functional.  The space $C(\Adel_\infty)$ is a free essential module
over~$\Adel_\infty$.  Since any smooth representation of~$\Adel_\fin$
is flat, $C(\Adel_\fin)$ is flat as a smooth representation
of~$\Adel_\GF$.  Thus we have constructed a resolution of the trivial
representation~$\C(1)$ in the category of smooth representations of
$\Adel\inv_\GF\cross\Adel_\GF$ that is flat as a representation
of~$\Adel_\GF$.  Hence we get an isomorphism
$$
\CHo(\Adel_\GF,\C(1)) \cong
\C(1)\Lhot_{\CCINF(\Adel_\GF)} \C(1) \cong
\C(1)\hot_{\CCINF(\Adel_\GF)} (C(\Adel_\infty),b) \cong
(\Lambda(\Adel_\infty),0)
$$
in the category $\Der_{\Adel\inv_\GF}$.  Notice that the boundary map
in this complex vanishes.

Putting everything together, we obtain a canonical isomorphism
$$
\CHH(\GF\inv\cross\Sch(\Adel_\GF))_1 \cong
\CHo(\GF\inv,\CHo(\Adel_\GF,\Sch(\Adel_\GF)_1)) \cong
\CHo(\GF\inv,\Sch(\Adel_\GF) \hot \Lambda \Adel_\infty).
$$

We now recall some results of~\cite{Meyer:Primes_Rep} that we are
going to need.  In~\cite{Meyer:Primes_Rep}, the
representation~$\lambda$ is used instead of~$\rho$.  This has no
effect because the Fourier transform intertwines the two
representations.  Let~$S$ be a sufficiently large finite set of
places.  Let $\ICL_S\defeq\Adel\inv_S/\GF\inv_S$ be the $S$\nbd{}idele
class group.  Let~$\Fourier$ be the Fourier transform on
$L^2(\Adel_S,dx)$.  It descends to a unitary operator on the Hilbert
space $L^2(\ICL_S,\abs{x}\,d\inv x)$.  Let
\begin{multline*}
  \Twist(\ICL_S) \defeq
  \{f\in L^2(\ICL_S,\abs{x}\,d\inv x) \mid \\
  \text{$f\abs{x}^\alpha\in\Sch(\ICL_S)$ and
    $\Fourier(f)\abs{x}^\alpha\in\Sch(\ICL_S)$ for all $\alpha>0$
  }\}.
\end{multline*}
It is shown in~\cite{Meyer:Primes_Rep} that
$$
\CHo(\GF\inv,\Sch(\Adel_\GF)) \cong
\Sch(\Adel_\GF)/\GF\inv \cong \varinjlim \Twist(\ICL_S)
$$
with suitable bornological embeddings
$\Twist(\ICL_S)\to\Twist(\ICL_{S'})$ for $S\subseteq S'$.  Moreover,
the right hand side can be identified with a space of functions
on~$\ICL_\GF$.  The following lemma allows us to apply this to
differential forms as well.

\begin{lemma}  \label{lem:trivialize_forms}
  There exists a bornological isomorphism on
  $\Sch(\Adel_\GF)\hot\Lambda \Adel_\infty$ that intertwines the
  representations $\rho\hot\pi$ and $\rho\hot\ID$.
\end{lemma}

\begin{proof}
  Let~$S$ be the set of infinite places of~$\GF$.
  We write
  $$
  \Sch(\Adel_\GF)\hot \Lambda\Adel_\infty=
  \Sch(\Adel_\fin) \hot \bigotimes_{v\in S} (\Sch(\GF_v) \hot
  \Lambda\GF_v).
  $$
  The representation of $\Adel\inv_\infty=\prod_{v\in S} \GF\inv_v$
  is the external tensor product representation.  Hence it suffices to
  construct such isomorphisms for Archimedean local fields, that
  is, for $\R$ and~$\C$.
  
  Consider first the case $\GF_v=\R$.  The spaces $\Lambda^0\R$ and
  $\Lambda^1\R$ are both $1$\nbd{}dimensional.  The action of~$\R\inv$
  on $\Lambda^0\R$ is trivial, on $\Lambda^1\R$ it is given by the
  identical quasi-character $x\mapsto x$.  Hence
  $\Sch(\R)\otimes\Lambda^1\R \cong \Sch(\R)\cdot x
  \subseteq\Sch(\R)$.  Decompose functions in $\Sch(\R)$ into an even
  and odd component $f=f_+ + f_-$ and let $\Phi(f)\defeq x df_+/dx +
  f_-$.  It is easy to see that this map is an isomorphism
  $\Sch(\R)\cong\Sch(\R)\cdot x$.
  
  Next we consider the case $\GF_v=\C$.  Again there is nothing to do
  for $\Lambda^0\C$.  The space $\Lambda^1\C$ has two generators
  corresponding to the holomorphic and antiholomorphic differentials
  $\partial/\partial z$ and $\partial/\partial \conj{z}$.  Since
  multiplication preserves the complex structure, they span invariant
  subspaces of $\Lambda^1\C$.  On the first subspace, $\C\inv$ acts by
  the identical quasi-character $z\mapsto z$, on the second by
  $z\mapsto\conj{z}$.  The space $\Lambda^2\C$ is again
  $1$\nbd{}dimensional and acted upon by the quasi-character $z\mapsto
  z\conj{z}$.  Hence we obtain the desired isomorphism
  for~$\Lambda^2\C$ by composing the isomorphisms for the two summands
  of~$\Lambda^1\C$.  For symmetry reasons, it suffices to treat
  $\Sch(\C)\partial/\partial y$.  We have to construct a
  $\C\inv$\nbd{}equivariant isomorphism $\Sch(\C)\cong\Sch(\C)\cdot
  z$.  We split $\Sch(\C)$ into homogeneous components with respect to
  the representation of the compact group $S^1\subseteq\C\inv$ and
  group together the representations $\conj{z}^n$ for $n\ge0$ into the
  positive and $z^n$ for $n>1$ into the negative part.  Thus we write
  any $f\in\Sch(\C)$ uniquely as $f=f_+ + f_-$.  The operator $\Phi
  f\defeq z\partial f_+/\partial z + f_-$.  is the desired
  isomorphism.
\end{proof}

It follows that
\begin{multline*}
  \CHH(\GF\inv\cross\Sch(\Adel_\GF))_1 \cong
  \CHo(\GF\inv,\Sch(\Adel_\GF) \hot \Lambda \Adel_\infty) \\ \cong
  (\Sch(\Adel_\GF) \hot \Lambda \Adel_\infty)/\GF\inv \cong
  \Sch(\Adel_\GF)/\GF\inv \hot \Lambda \Adel_\infty.
\end{multline*}
The expressions in the second line are graded bornological vector
spaces viewed as chain complexes with vanishing boundary map.  Hence
their graded components are the homogeneous parts of the Hochschild
homology.  We also get the additional fact that the homology
computation can be done using only homotopy equivalences.  Therefore,
it immediately implies results about Hochschild cohomology and also
about Hochschild homology and cohomology for tensor products with
$\GF\inv\cross\Sch(\Adel_\GF)$.  The last isomorphism is not quite
natural but can be chosen $\Adel\inv_\GF$\nbd{}equivariant.  We can
identify $\Sch(\Adel_\GF)/\GF\inv$ further with $\varinjlim
\Twist(\ICL_S)$ as in~\cite{Meyer:Primes_Rep}.

\subsection{The inhomogeneous part of the cyclic homology theories}
\label{sec:explicit_inhomogeneous}

In order to pass to cyclic and periodic cyclic homology, we need
explicit chain maps that implement the above homology computation.
First we consider the inhomogeneous part.  The integration map
$(\Sch(\Adel_\GF),*)\to\C$ is an algebra homomorphism.  Its Fourier
transform is the map $(\Sch(\Adel_\GF),\cdot)\to\C$ of evaluation
at~$0$.  We get an $\Adel\inv_\GF$\nbd{}equivariant algebra
homomorphism
$$
\Phi\colon \GF\inv\cross\Sch(\Adel_\GF)\to
\GF\inv\cross\C(1) = \C[\GF\inv].
$$

\begin{proposition}  \label{pro:inhom}
  The map~$\Phi$ induces a homotopy equivalence
  $$
  \CHH(\GF\inv\cross\Sch(\Adel_\GF))_\inhomogen \to
  \CHH(\C[\GF\inv])_\inhomogen
  $$
  and hence isomorphisms between the inhomogeneous parts of the
  Hochschild homology, the cyclic homology and the periodic cyclic
  homology of $\GF\inv\cross\Sch(\Adel_\GF)$ and $\C[\GF\inv]$.  The
  group~$\ICL_\GF$ acts trivially on these homology spaces.
\end{proposition}

This is essentially proved by our above computations.  We make the
assertion more explicit by describing the Hochschild homology of
$\C[\GF\inv]$.  Let $\Roots\subseteq\GF\inv$ be the subgroup of roots
of unity.  The group $\GF\inv/\Roots$ is a free Abelian group of
infinite rank.  Hence we have an isomorphism
$\GF\inv\cong\Roots\times\Z^\N$.  Therefore,
$$
\C[\GF\inv]\cong\C[\Roots]\otimes \C[\Z]^{\otimes\N},
$$
where we define the infinite tensor product as an inductive limit of
finite tensor products.  Each factor $\C[\Z]$ is a Laurent series
ring.  We view $\C[\GF\inv]$ as the algebra of Laurent series on the
Pontrjagin dual $\widehat{\GF\inv}$.  We let
$\Omega(\widehat{\GF\inv})$ be the algebra of differential forms
on~$\widehat{\GF\inv}$ with Laurent series as coefficients.  This is a
differential graded algebra over $\C[\GF\inv]$, so that we obtain a
canonical homomorphism of differential graded algebras
$\Omega(\C[\GF\inv])\to\Omega(\widehat{\GF\inv})$.  It is well-known
that the map
\begin{equation}  \label{eq:inhom_explicit}
  \CHH(\C[\GF\inv]) =
  (\Omega(\C[\GF\inv]),b)\to
  (\Omega(\widehat{\GF\inv}),0)
\end{equation}
is a homotopy equivalence.  Hence the Hochschild homology of
$\C[\GF\inv]$ is $\Omega(\widehat{\GF\inv})$.

For any $g\in\GF\inv$, the $g$\nbd{}homogeneous subspace of
$\Omega^n(\widehat{\GF\inv})$ is canonically isomorphic to the group
homology of~$\GF\inv$.  When we computed the inhomogeneous part of
$\CHH(\GF\inv\cross\Sch(\Adel_\GF))$, the only important ingredient
was exactly the integration map $\Sch(\Adel_\GF)\to\C$.  This shows
that the isomorphism from $\CHH(\GF\inv\cross\Sch(\Adel_\GF))_g$ to
$\CHo(\GF\inv,\C(1))$ that we constructed above agrees with the map
induced by~$\Phi$.  Since~$\Phi$ is an algebra homomorphism, it
induces maps on cyclic and periodic cyclic homology.  These maps are
compatible with the decomposition into homogeneous and inhomogeneous
parts.  Since~$\Phi$ induces an isomorphism on the inhomogeneous part
of the Hochschild homology, it also induces isomorphisms on the
inhomogeneous parts of the cyclic and periodic cyclic homologies.
Since~$\Phi$ is invariant under~$\Adel\inv_\GF$, the action
of~$\ICL_\GF$ on the $g$\nbd{}homogeneous part of the Hochschild
homology is trivial.  This finishes the proof of
Proposition~\ref{pro:inhom}.

Let $N$ and~$d_\dR$ be the number operator (that is, $N=n$ on
$n$\nbd{}forms) and the de Rham differential on
$\Omega(\widehat{\GF\inv})$.  Then the homotopy equivalence
in~\eqref{eq:inhom_explicit} intertwines $B$ and $N d_\dR$.  Hence we
obtain a chain map from the $(B,b)$\nbd{}bicomplex to the bicomplex
$(\Omega(\widehat{\GF\inv}),N d_\dR,0)$ that is a homotopy equivalence
with respect to the vertical boundary maps $b$ and~$0$.  Therefore, it
is a homotopy equivalence between the total complexes as well.  Since
the total complex of the $(B,b)$\nbd{}bicomplex computes the cyclic
homology of $\C[\GF\inv]$, we get
$$
\HC_j(\C[\GF\inv]) \cong
\Omega^j(\widehat{\GF\inv})/d_\dR\Omega^{j-1}(\widehat{\GF\inv})
\oplus
H^{j-2}_\dR(\widehat{\GF\inv})
\oplus
H^{j-4}_\dR(\widehat{\GF\inv})
\oplus \cdots
$$
for all $j\in\N$, where $H^j_\dR$ denotes the de Rham cohomology
of~$\widehat{\GF\inv}$.  The periodicity operator~$S$ can also be
described easily, so that we get
$$
\HP_j(\C[\GF\inv]) = \prod_{k\in\Z}
H_\dR^{j+2k}(\widehat{\GF\inv}).
$$
We leave it to the reader to describe the coaction and decomposition
into homogeneous and inhomogeneous parts of the Hochschild and cyclic
homology.  The coaction of~$\Z^\N$ on $\HP_j(\C[\GF\inv])$ is trivial.
However, we still have a non-trivial coaction of~$\Roots$.  An element
of $H_\dR^j(\widehat{\GF\inv})$ is homogeneous if and only if it is
invariant under the action of~$\widehat{\Roots}$ that permutes the
connected components of~$\widehat{\GF\inv}$.

\subsection{The homogeneous part of the cyclic homology theories}
\label{sec:hom_HH}

We can also describe an explicit chain map implementing our
computation of the homogeneous part of the Hochschild homology.  It is
more convenient to equip $\Sch(\Adel_\GF)$ with the pointwise product
in order to bring into play the de Rham boundary map.

Let $\Sch\Omega(\Adel_\infty)$ be the differential graded algebra of
differential forms on $\Adel_\infty\cong\R^{n_\GF}$ with Schwartz
functions as coefficients.  Let $\Sch\Omega(\Adel_\GF)\defeq
\Sch\Omega(\Adel_\infty)\hot\Sch(\Adel_\fin)$.  We let $N$
and~$d_\dR$ be the number operator and the de Rham boundary maps on
these spaces of differential forms.  We let~$\Adel\inv_\GF$ act on
$\Sch\Omega(\Adel_\GF)$ in the canonical way.  Thus we can form the
coinvariant space $\Sch\Omega^n(\Adel_\GF)/\GF\inv$.  The Fourier
transform gives rise to a natural isomorphism
$$
\Fourier \colon \Sch\Omega(\Adel_\GF) \congto
\Sch(\Adel_\GF) \hot \Lambda \Adel_\infty.
$$
Hence our computation of the homogeneous part of the Hochschild
homology yields
$$
\HH_n(\GF\inv\cross\Sch(\Adel_\GF))_1 \cong
(\Sch(\Adel_\GF) \hot \Lambda \Adel_\infty)/\GF\inv \cong
\Sch\Omega(\Adel_\GF)/\GF\inv.
$$
We can realize this by the following explicit map.  The homogeneous
part of $\Omega^n(\GF\inv\cross\Sch(\Adel_\GF))$ is the closed linear
span of monomials
$$
f_0\delta_{g_0}\,d(f_1\delta_{g_1})\dots d(f_n\delta_{g_n}),
\qquad
d(f_1\delta_{g_1})\dots d(f_n\delta_{g_n}),
$$
with $g_0\cdots g_n=1$ and $g_1\cdots g_n=1$, respectively.  We map
them to the classes of the differential forms
$$
f_0 d(\lambda_{g_0} f_1)\wedge \dots \wedge
d(\lambda_{g_0\cdots g_{n-1}} f_n),
\qquad
df_1\wedge d(\lambda_{g_1} f_2) \dots \wedge
d(\lambda_{g_1\cdots g_{n-1}} f_n)
$$
in $\Sch\Omega^n(\Adel_\GF)/\GF\inv$, respectively.  The resulting map
$$
\Phi_1\colon
\Omega(\GF\inv\cross\Sch(\Adel_\GF))_1 \to
\Sch\Omega(\Adel_\GF)/\GF\inv
$$
evidently satisfies $\Phi_1\circ b=0$ and $\Phi_1\circ B=Nd_\dR\circ
\Phi_1$.  A closer inspection of our above computations shows
that~$\Phi_1$ is a homotopy equivalence
$$
(\Omega(\GF\inv\cross\Sch(\Adel_\GF)),b)_1 \to
(\Sch\Omega(\Adel_\GF)/\GF\inv,0).
$$
Hence we can use~$\Phi_1$ also to compute the cyclic and periodic
cyclic homology spaces.  As in the computation for the inhomogeneous
part, we get
$$
\HP_j(\GF\inv\cross\Sch(\Adel_\GF))_1 \cong
\bigoplus_{k\in\Z}
H_{j-2k}(\Sch\Omega(\Adel_\GF)/\GF\inv,(d_\dR)_\ast)
$$
and a similar result for the cyclic theories.

It remains to compute this homology more explicitly.  We use once
again that $\CHo(\GF\inv,\Sch\Omega^n(\Adel_\GF))\cong0$.  This means
that the complex $(\Sch\Omega(\Adel_\GF),d_\dR)$ is acyclic with
respect to the group homology functor.  We want to describe its
homology.  Fix an orientation on~$\Adel_\infty$ and let
$$
{\smallint}\colon
\Sch\Omega^{n_\GF}(\Adel_\GF) \cong
\Sch\Omega(\Adel_\infty)\hot\Sch(\Adel_\fin) \to
\Sch(\Adel_\fin)
$$
be the usual integration map for forms of top degree.  Let
$\sigma\colon\Adel\inv_\GF\to\{\pm1\}$ be the character that
multiplies the signs of an idele at the real places of~$\GF$.  Thus
$\sigma(a)=+1$ if multiplication by~$a$ preserves orientation and
$\sigma(a)=-1$ if multiplication by~$a$ reverses orientation.
Hence~${\smallint}$ is an equivariant map to
$\Sch(\Adel_\fin)\hot\C(\sigma)$.  We claim that
$(\Sch\Omega(\Adel_\GF),d_\dR,{\smallint})$ is a linearly split
resolution of $\Sch(\Adel_\fin)\hot\C(\sigma)$.  To prove this we
write
$$
(\Sch\Omega(\Adel_\GF),d_\dR) \cong
(\Sch(\Adel_\fin),0) \hot (\Sch\Omega(\R),d_\dR)^{\hot n_\GF}.
$$
It is elementary to verify that $\Sch(\R)
\overset{d_\dR}\longrightarrow \Sch\Omega^1(\R)
\overset{\smallint}\longrightarrow \C$ is a linearly split extension.
This implies the claim.  Since the complex
$(\Sch\Omega(\Adel_\GF),d_\dR)$ is also acyclic for the group homology
functor, we obtain
$$
H_j(\Sch\Omega(\Adel_\GF)/\GF\inv,(d_\dR)_\ast) \cong
H_{n_\GF-j}(\GF\inv,\Sch(\Adel_\fin)\hot\C(\sigma)).
$$

Now we use the $\Adel\inv_\fin$\nbd{}equivariant isomorphism
$\Sch(\Adel_\fin)\cong\CCINF(\Adel\inv_\fin)$ constructed
in~\cite{Meyer:Primes_Rep}.  Let $(U_k)_{k\in\N}$ be a decreasing
sequence of compact neighborhoods of the identity in~$\Adel\inv_\fin$
and let
$$
L_k\defeq \GF\inv\cap U_k\subseteq \GF\inv.
$$
Since $g\in L_k$ implies $\abs{g}_v=1$ for all finite places~$v$, the
groups~$L_k$ are contained in~$\GF\inv_\infty$.  One can show that
they form a decreasing sequence of cocompact subgroups
of~$\GF\inv_\infty$ with $\bigcap L_k=\{1\}$.  Choosing~$U_1$
sufficiently small, we can achieve that~$H_1$ contains no roots of
unity and that $\sigma|_{L_1}=1$.  It follows that~$L_k$ is a free
Abelian group of rank~$u$ for all~$n$, where~$u$ is the rank
of~$\GF\inv_\infty$.

The space $\CCINF(\Adel\inv_\fin)$ is the direct union of the
$\Adel\inv_\fin$\brd{}invariant subspaces
$\CCINF(\Adel\inv_\fin/U_k)$.  Hence
$$
H_*(\GF\inv,\CCINF(\Adel\inv_\fin)\otimes\C(\sigma)) \cong
\varinjlim H_*(\GF\inv,\CCINF(\Adel\inv_\fin/U_k)\otimes\C(\sigma)).
$$
For any $k\in\N$, the quotient space $\Adel\inv_\fin/U_k$ is
discrete and~$L_k$ acts trivially on it.  The induced action of
$\GF\inv/L_k$ on~$\Adel\inv_\fin/U_k$ is free.  Hence the Shapiro
Lemma yields
$$
H_*(\GF\inv,\CCINF(\Adel\inv_\fin/U_k)\otimes\C(\sigma)) \cong
H_*(L_k) \otimes
(\CCINF(\Adel\inv_\fin/U_k)\otimes\C(\sigma))/\GF\inv.
$$
The latter tensor factor is naturally isomorphic to the space
$\CCINF(\Adel\inv_\fin/U_k\GF\inv)_\sigma$ of smooth functions
$f\colon\Adel\inv_\fin/U_k\to\C$ that satisfy $f(ax)=\sigma(a)f(x)$
for all $a\in\GF\inv$, $x\in\Adel\inv_\fin/U_k$.  We can identify it
unnaturally with $\CCINF(\Adel\inv_\fin/U_k\GF\inv)$ by choosing a
fundamental domain for the action of~$\GF\inv$.  The group homology of
the free Abelian group~$L_k$ is an exterior algebra on~$u$ generators.

Inspection shows that the passage from~$U_k$ to~$U_{k+1}$ is the
tensor product of the obvious map on $\CCINF(\Adel\inv_\fin/\GF\inv
U_k)_\sigma$ and multiplication with certain constants that depend on
the group $L_k/L_{k+1}$ on $H_*(L_k)$.  Thus
$$
H_*(\GF\inv,\CCINF(\Adel\inv_\fin)\otimes\C(\sigma)) \cong
\varinjlim \Lambda^\ast \C^u \otimes \CCINF(\Adel\inv_\fin/\GF\inv
U_k)_\sigma
\cong \Lambda^\ast \C^u \otimes
\CCINF(\Adel\inv_\fin/\cl{\GF\inv})_\sigma.
$$
Here~$\cl{\GF\inv}$ denotes the closure of~$\GF\inv$
in~$\Adel\inv_\fin$.  That is, we obtain the space of functions
$f\colon \Adel\inv_\fin\to\Lambda^\ast \C^u$ that satisfy
$f(ax)=\sigma(a) f(x)$ for all $a\in\GF\inv$, $x\in\Adel\inv_\fin$.
The action of~$\Adel\inv_\fin$ on the right hand side is trivial on
the first and the obvious action on the second tensor factor.  The
group~$\Adel\inv_\infty$ acts by~$\sigma$.

This finishes our computation of the cyclic and periodic cyclic
homology of the crossed product.

\section{The topological K-theory of the crossed product}
\label{sec:Ktop}

Let~$A$ be the (complex) $C^*$\nbd{}algebraic crossed product
$\GF\inv\cross C_0(\Adel_\GF)$.  We want to compute the topological
$\K$\nbd{}theory $\K(A)$.  The group~$\GF\inv$ is Abelian and hence
amenable.  Therefore, the full and reduced crossed products agree
and~$\GF\inv$ satisfies the Baum-Connes conjecture with arbitrary
coefficients by a result of Nigel Higson and Gennadi Kasparov
(\cite{Higson-Kasparov:Amenable}).  Thus we can replace $\K(A)$ by
$\Ktop(\GF\inv,C_0(\Adel_\GF))$.  The latter group can always be
computed by a spectral sequence.  However, in order to get more
precise results it is often better to use special structure present in
the situation.

Let $\Roots\subseteq\GF\inv$ be the subgroup of roots of unity and
choose an isomorphism $\GF\inv\cong\Roots\times\Z^\N$ as above.  Let
$$
B\defeq \Roots\cross C_0(\Adel_\GF),
$$
then we have $A\cong \Z^\N\cross B$ and hence $\K(A)\cong
\Ktop(\Z^\N,B)$.  Since~$\Z^\N$ is torsion free, the spectral sequence
for $\Ktop(\Z^\N,B)$  is much simpler than for $\Ktop(\GF\inv,A)$.
Its $E^2$\nbd{}term is
$$
E^2_{pq} \defeq H_p(\Z^\N, \K_q(B)).
$$
Thus our first step is to compute $\K(\Roots\cross C_0(\Adel_\GF))
\cong \K^\Roots(C_0(\Adel_\GF))$ together with the representation of
$\GF\inv/\Roots$ on it.  To express the result in a uniform way for
all global fields, we need some ingredients.  They all are
$\Ztwo$\nbd{}graded Abelian groups together with a representation
of~$\Adel\inv_\GF$.

Let $\Rep(\Roots)$ be the representation ring of~$\Roots$ in even
degree and~$0$ in odd degree.  Thus $\Rep(\Roots)\cong
\K^\Roots(\C)$.  The forgetful map $\K^\Roots(\C)\to\K(\C)$ associates
to each representation its dimension.  We let
$J\Rep(\Roots)\subseteq\Rep(\Roots)$ be its kernel.  We
let~$\Adel\inv_\GF$ act trivially on $\Rep(\Roots)$ and
$J\Rep(\Roots)$.  Let $n=n_\GF$ be the degree of~$\GF$ over~$\Q$, so
that $\Adel_\infty\cong\R^n$.  Let~$\sigma$ be the orientation
character that already occured in the periodic cyclic homology and let
$\Z(\sigma)$ be the resulting representation of~$\Adel\inv_\GF$
on~$\Z$.  We view $\Z(\sigma)$ as a $\Ztwo$\nbd{}graded Abelian group
in degree $n \bmod (2)$.  For a totally disconnected locally compact
space~$F$ and a graded Abelian group~$V$, we let $\CCINF(F,V)$ be the
space of locally constant, compactly supported functions $F\to V$ with
the obvious grading.  If a group acts on~$F$ and~$V$, we equip
$\CCINF(F,V)$ with the induced action.  In particular, on
$\CCINF(\Adel_\fin,\Z(\sigma))$, the finite adeles act by translation
and the infinite adeles act by~$\sigma$.  Finally, we write~$V^\Roots$
for the subspace of $\Roots$\nbd{}invariants in a representation~$V$
of~$\Roots$.

\begin{proposition}  \label{pro:K_Roots_Adel}
  There is an $\Adel\inv_\GF$\brd{}equivariant grading preserving
  isomorphism 
  $$
  \K^\Roots(C_0(\Adel_\GF)) \cong
  J\Rep(\Roots) \oplus
  \CCINF(\Adel_\fin,\Z(\sigma))^\Roots.
  $$
  Evaluation at $0\in\Adel_\GF$ induces a map
  $\K^\Roots(C_0(\Adel_\GF))\to\K^\Roots(\C)$, which sends the first
  summand isomorphically onto $J\Rep(\Roots)\subseteq\K^\Roots(\C)$.
\end{proposition}

\begin{proof}
  Let $\Adel^*_\fin\defeq \Adel_\fin\setminus\{0\}$.  Evaluation at
  $0\in\Adel_\fin$ gives rise to an extension of $C^*$\nbd{}algebras
  $$
  C_0(\Adel^*_\fin \times \Adel_\infty) \into
  C_0(\Adel_\GF) \prto
  C_0(\Adel_\infty).
  $$
  This extension splits by the following $\Roots$\nbd{}equivariant
  $*$\nbd{}homomorphism.  Let $\MCR\subseteq\Adel_\fin$ be the
  maximal compact subring.  Thus $x\in\MCR$ if and only if
  $\abs{x}_v\le1$ for all finite places~$v$.  Let~$1_\MCR$ be the
  characteristic function of~$\MCR$.  The map $f\mapsto
  1_\MCR\otimes f$ is the desired section.  Applying the
  $\K$\nbd{}theory functor, we get a split exact sequence
  $$
  \K^\Roots(C_0(\Adel^*_\fin \times \Adel_\infty)) \into
  \K^\Roots(C_0(\Adel_\GF)) \prto
  \K^\Roots(C_0(\Adel_\infty)).
  $$
  Since $\Adel_\infty\cong\R^n$, we would like to use equivariant Bott
  periodicity to simplify the above expressions.  For that we need to
  know whether~$\R^n$ has a $\Roots$\nbd{}invariant spin structure.
  We distinguish three cases.

  Suppose first that~$\GF$ has no real places.  Thus $\sigma=1$
  and~$n$ is even.  The space~$\Adel_\infty$ is a complex manifold in
  a natural and in particular $\Adel\inv_\infty$\nbd{}invariant way.
  Since a complex structure is more than a spin structure, equivariant
  Boot periodicity applies and shows that $C_0(\Adel_\infty)$ is
  $\KK^\Roots$\nbd{}equivalent to~$\C$.  Therefore,
  \begin{align*}
    \K^\Roots(C_0(\Adel_\infty))
    &\cong \K^\Roots(\C)
    \cong \Rep(\Roots),
    \\
    \K^\Roots(C_0(\Adel^*_\fin \times \Adel_\infty))
    &\cong \K^\Roots(C_0(\Adel^*_\fin))
    \cong \K(C_0(\Adel^*_\fin/\Roots))
    \\
    &\cong \CCINF(\Adel^*_\fin/\Roots,\Z)
    \cong \CCINF(\Adel^*_\fin,\Z)^\Roots.
  \end{align*}
  We have also used that~$\Roots$ acts freely on $\Adel^*_\fin$ and
  that $\Adel^*_\fin/\Roots$ is totally disconnected.  Since the
  action of~$\Adel\inv_\infty$ preserves the spin structure
  on~$\Adel_\infty$, the above isomorphisms are
  $\Adel\inv_\GF$\nbd{}equivariant.  They also preserve the grading.
  
  Although the $C^*$\nbd{}algebra extension above splits
  $\Roots$\nbd{}equivariantly, the section we constructed is
  manifestly not equivariant for the action of~$\Adel\inv_\GF$.
  Nevertheless, the restriction of the induced map on $\K$\nbd{}theory
  to $J\Rep(\Roots)$ is equivariant.  The quotient
  $\Rep(\Roots)/J\Rep(\Roots)\cong\Z$ has the effect of putting back
  the point $0\in\Adel_\fin$.  Thus we obtain the desired description
  of $\K^\Roots(C_0(\Adel_\GF))$.  It is left to the reader to check
  that $J\Rep(\Roots)$ is detected by evaluation at~$0$.

  Next suppose that~$\GF$ has real places and that~$n$ is even.
  Equivalently, the number of real places is even.  Then
  $\Roots=\{\pm1\}$.  The action of~$-1$ on $\R^2$ is the same as the
  action of~$-1$ on~$\C$ and hence again preserves a complex
  structure.  Hence we can apply equivariant Bott periodicity as in
  the purely complex case and get essentially the same results.  The
  only difference is that~$\Adel\inv_\infty$ need not preserve the
  spin structure on~$\Adel_\infty$ and therefore acts by the
  orientation character~$\sigma$.  This yields the assertions if~$n$
  is even.

  Finally, suppose that~$n$ is odd.  Hence there is an odd number of
  real places.  Again this forces $\Roots=\{\pm1\}$.  Equivariant Bott
  periodicity allows us, as above, to remove all complex places and an
  even number of real places from~$\Adel_\infty$.  Thus
  $C_0(\Adel_\infty)$ is $\KK^\Roots$\brd{}equivalent to $C_0(\R)$
  with the action of~$\Roots$ by reflection at the origin.  A
  straightforward computation shows that evaluation at~$0$ induces an
  isomorphism $\K^\Roots(C_0(\R))\cong
  J\Rep(\Roots)\subseteq\Rep(\Roots)$.  There exists a fundamental
  domain~$F$ for the free action of~$\Roots$ on the totally
  disconnected space~$\Adel^*_\fin$.  This allows us to identify
  $$
  \K^\Roots(C_0(\Adel_\fin^*\times\Adel_\infty)) \cong
  \K^\Roots(C_0(\Roots\times F\times\Adel_\infty)) \cong
  \K(F\times\Adel_\infty) \cong
  \CCINF(F,\Z(\sigma)).
  $$
  Recall that $\CCINF(F,\Z(\sigma))$ sits in odd degree in this
  case.  As above, the infinite adeles act on this by~$\sigma$.  Since
  $\sigma(-1)=-1$, we should extend functions on~$F$ to~$\Adel^*_\fin$
  by $f(-x)=-f(x)$.  Since this condition forces $f(0)=0$, we are
  allowed to put the point~$0$ back and get the required result also
  in this case.
\end{proof}

Define $\Ztwo$\nbd{}graded groups by
$$
\K_\homogen \defeq
(\CCINF(\Adel_\fin,\Z)\otimes\Z(\sigma))^\Roots,
\qquad
\K_\inhomogen \defeq
J\Rep(\Roots).
$$
They correspond to the homogeneous and inhomogeneous parts in the
cyclic homology computations.  We obtain
$$
E^2_{pq} \cong H_p(\GF\inv/\Roots,\K_\homogen) \oplus
H_p(\GF\inv/\Roots,\K_\inhomogen).
$$

\begin{lemma}  \label{lem:Ktop_collapse}
  The boundary maps $d^2_{pq}\colon E^2_{pq}\to E^2_{p-2,q+1}$ vanish
  for all $p,q$.
\end{lemma}

\begin{proof}
  If~$n$ is even, this is clear because either the source or target
  of~$d^2_{pq}$ vanish.  Let~$n$ be odd and let~$g$ be the adele
  $(1_\fin,-1_\infty)$.  That is, the value of~$g$ at a place~$v$
  is~$\pm1$ depending on whether~$v$ is infinite or not.  Then~$g$
  acts as~$+1$ on $J\Rep(\Roots)$ and as~$(-1)^n=-1$ on the other
  summand.  Since the boundary map~$d^2_{pq}$ is natural with respect
  to equivariant $*$\nbd{}homomorphisms, it commutes with the action
  of~$g$.  Hence it cannot mix the two direct summands.  Since each
  summand is concentrated in one parity, we get $d^2_{pq}=0$.
\end{proof}

Lemma~\ref{lem:Ktop_collapse} implies that there is a filtration on
$\K_m(\GF\inv\cross C_0(\Adel_\GF))$ whose subquotients are isomorphic
to~$E^2_{pq}$ for $p+q=m$.  Since~$\GF\inv$ acts trivially on
$\K_\inhomogen$, we get natural isomorphisms
$$
H_p(\GF\inv/\Roots,\K_\inhomogen) \cong
\K_\inhomogen \otimes H_p(\Z^\N) \cong
J\Rep(\Roots) \otimes \Lambda^p \Z^\N,
$$
where $\Lambda^p V$ is the $p$th exterior power of~$V$.  The
evaluation map at~$0$ maps this group isomorphically onto the
corresonding subspace of $\K(\GF\inv\cross\C)\cong
\K(C_0(\widehat{\GF\inv}))$.

The contribution of $\K_\homogen$ requires considerably more work.  We
proceed in two steps.  We first deal with the subgroup
$\GF\inv_\infty/\Roots$ and then use the spectral sequence 
$$
H_p(\GF\inv/\GF\inv_\infty,H_q(\GF\inv_\infty/\Roots, \K_\homogen))
\Longrightarrow H_{p+q}(\GF\inv/\Roots, \K_\homogen).
$$

Let $\MCR\subseteq\Adel_\fin$ be the maximal compact subring, that is,
the set of all $x\in\Adel_\fin$ with $\abs{x}_v\le1$ for all finite
places~$v$.  There is a basis for the neighborhoods of zero
in~$\Adel_\fin$ consisting of compact, open subgroups of the form
$x_m\cdot\MCR$ for a suitable sequence of ideles~$(x_m)$.  Thus
$$
\K_\homogen = \CCINF(\Adel_\fin,\Z(\sigma))^\Roots \cong
\varinjlim \CCINF(\Adel_\fin/x_m\MCR,\Z(\sigma))^\Roots.
$$
The action of $\GF\inv_\infty/\Roots$ respects this inductive limit
decomposition, so that
$$
H_*(\GF\inv_\infty/\Roots, \K_\homogen)
\cong \varinjlim H_*(\GF\inv_\infty/\Roots,
\CCINF(\Adel_\fin/x_m\MCR,\Z(\sigma))^\Roots).
$$

The space $\Adel_\fin/x_m\MCR$ is discrete.  We let
$$
L_m(y)\defeq \{g\in\GF\inv_\infty\mid gy-y\in x_m\MCR\}
$$
be the stabilizer of the point $y+x_m\MCR\in \Adel_\fin/x_m\MCR$.
This is a cocompact subgroups of $\GF\inv_\infty$ because the action
of~$\GF\inv_\infty$ on~$\Adel_\fin$ factors through the compact group
$\cl{\GF\inv_\infty}\subseteq\MCR\inv\subseteq \Adel\inv_\fin$.  Let
$\dot{L}_m(y)\defeq L_m(y)/(\Roots\cap L_m(y))$.  This is a cocompact
subgroup of $\GF\inv_\infty/\Roots$ and hence a free Abelian group of
rank~$u$.

The space $\Adel_\fin/x_m\MCR$ is a disjoint union of orbits
$\GF\inv_\infty/L_m(y)$.  Each orbit yields a contribution
$$
H_*(\GF\inv_\infty/\Roots,
\CCINF(\GF\inv_\infty/L_m(y),\Z(\sigma))^\Roots).
$$
To compute the latter, we have to distinguish three cases.

If $-1\in L_m(y)$ and $\sigma(-1)=-1$, then there simply are no
$\Roots$\nbd{}invariant functions
$\GF\inv_\infty/L_m(y)\to\Z(\sigma)$, so that we get~$0$.  In the
other cases, $\sigma$ vanishes on $\Roots\cap L_m(y)$.  Hence
$\CCINF(\GF\inv_\infty/L_m(y),\Z(\sigma))^\Roots$ is isomorphic to
$\CCINF(\GF\inv_\infty/L_m(y)\Roots,\Z(\sigma))$.  This
representation is induced from the representation $\Z(\sigma)$ of
$\dot{L}_m(y)$.  Hence we obtain
$$
H_*(\GF\inv_\infty/\Roots,
\CCINF(\GF\inv_\infty/L_m(y)\Roots,\Z(\sigma))) \cong
H_*(\dot{L}_m(y),\Z(\sigma)).
$$
If~$\sigma$ does not vanish on $L_m(y)$, we obtain again~$0$.
Otherwise, we obtain the group homology of $\dot{L}_m(y)$, which is
$\Lambda\dot{L}_m(y)$.  The final result is therefore
$$
H_*(\GF\inv_\infty/\Roots,
\CCINF(\Adel_\fin/x_m\MCR,\Z(\sigma))^\Roots) \cong
\bigoplus_{\{y\in\Adel_\fin/x_m\MCR\mid L_m(y)\subseteq\Ker\sigma\}}
\Lambda\dot{L}_m(y).
$$

To understand how these spaces fit together in the inductive limit
over~$m$, we have to rewrite this result in a more functorial way.
Let~$L$ be a free Abelian group of rank~$u$.  Let
$\omega\in\Hom(\Lambda^u L,\Z)\cong \Z$ be one of the two generators.
We define a biadditive pairing $\Lambda^p L\otimes \Lambda^{u-p}
L\to\Z$ by
$$
(\xi,\eta) \defeq \omega(\xi\wedge\nu),
$$
This pairing defines an isomorphism $\Lambda^p L\cong\Hom(\Lambda^{u-p}
L,\Z)$.

Let $L'\subseteq L$ be a cocompact subgroup of a free Abelian group.
Then we can consider the subspace $\Lambda(L':L)$ of all
$f\in\Lambda^p L\otimes\Q$ with $(f,x)\in\Z$ for all $x\in
\Lambda^{n-p} L'$.  This subspace does not depend on the choice
of~$\omega$.  Explicitly, it looks as follows.  Choose a basis
$(g_1,\dots,g_u)$ for~$L$ such that $(n_1g_1,\dots,n_ug_u)$ is a basis
for~$L'$ with certain $n_j\in\N\inv$.  Express all exterior forms in
these bases.  Then $\Lambda(L':L)$ is equal to the subalgebra
generated by the $1$\nbd{}forms $n_j^{-1} dg_j$.  This subspace is
non-canonically isomorphic to $\Lambda^p L'$.

For each $y\in\Adel_\fin$, we have a decreasing sequence of subgroups
$\dot{L}_m(y)\subseteq\GF\inv_\infty/\GF\inv$.  Let
$$
\Lambda_y\defeq
\bigcup_{m\in\N} \Lambda(\dot{L}_m(y):\GF\inv_\infty/\Roots)
\subseteq
\Lambda(\GF\inv_\infty/\Roots) \otimes\Q.
$$
We can now reformulate the above computation as follows:

\begin{proposition}  \label{pro:integer_GF_homology}
  The homology $H_*(\GF\inv_\infty/\Roots, \K_\homogen)$ is naturally
  isomorphic to the space of all functions in $\CCINF(\Adel_\fin,
  \Lambda(\GF\inv_\infty/\Roots)\otimes
  \Q(\sigma))^{\cl{\GF\inv_\infty}}$ that satisfy the integrality
  condition $f(y)\in\Lambda_y$ for all $y\in\Adel_\fin$.
\end{proposition}

\begin{proof}
  The space of functions that satisfies the integrality condition is
  the inductive limit of the spaces $\CCINF(\Adel_\fin/x_m\MCR,
  \Lambda(\dot{L}_m(y):\GF\inv_\infty)
  \otimes\Z(\sigma))^{\cl{\GF\inv_\infty}}$ for $m\to\infty$.  For a
  fixed~$m$ this space is isomorphic to the result of our computation
  for the group homology.  The point is that these isomorphisms can be
  chosen natural and compatible with the inductive limit.  Let
  $L\defeq \GF\inv_\infty/\Roots$.  The group homology of the free
  Abelian group~$L$ with coefficients~$\K_\homogen$ is computed by the
  Koszul complex $(\K_\homogen\otimes\Lambda(L),\delta)$ with
  $$
  \delta(f\,dg_1\wedge\dots\wedge dg_k) \defeq
  \sum_{j=1}^k (-1)^j (f-g_j\cdot f) \,dg_1\wedge\dots \wedge
  \widehat{dg_j}\wedge\dots\wedge dg_k.
  $$
  Consider the graded $\Q$\nbd{}vector space
  $$
  W\defeq  \CCINF(\Adel_\fin,
  \Lambda(L)\otimes\Q(\sigma))^{\cl{\GF\inv_\infty}}
  $$
  as a complex with trivial boundary map.
  
  The closure $\cl{\GF\inv_\infty}\subseteq\Adel\inv_\fin$
  of~$\GF\inv_\infty$ is a compact group.  Let~$\mu$ be the normalized
  Haar measure on $\cl{\GF\inv_\infty}/\Roots$.  We extend~$\sigma$ to
  a continuous character~$\cl{\sigma}$ on~$\cl{\GF\inv_\infty}$ and
  consider the integration map $\alpha\colon \K_\homogen\otimes
  \Lambda(L) \to W$ defined by
  $$
  \alpha (f\,dg_{i_1}\wedge\dots\wedge dg_{i_k})(t) \defeq
  \int_{\cl{\GF\inv_\infty}/\Roots} f(at) \cl{\sigma}(a) \,d\mu(a)
  \,dg_{i_1}\wedge\dots\wedge dg_{i_k}.
  $$
  It is a chain map because the integration annihilates all
  coinvariants $f-g\cdot f$ for $g\in\GF\inv_\infty$.  The map on
  homology induced by~$\alpha$ is injective.  This follows easily from
  our elementary computations with $\Adel_\fin/x_m\MCR$.  The
  assertion follows by explicitly identifying the image of
  $H_*(\K_\homogen\otimes\Lambda(L))$ in~$W$.
\end{proof}

Let~$W$ be the graded Abelian group that we identified with the
homology in Proposition~\ref{pro:integer_GF_homology}.  We still have
to describe the homology of $\GF\inv/\GF\inv_\infty$ with coefficients
in~$W$.  We use a filtration of~$W$ for this purpose.  Let~$S$ be a
finite set of finite places.  Let~$T$ be the set of finite places not
in~$S$, let~$\MCR_T$ be the set of all $x\in\Adel_T$ that satisfy
$\abs{x}_v\le1$ for all $v\in T$ and let $1_{\MCR_T}\in \Sch(\Adel_T)$
be its characteristic function.

A smooth function on~$\Adel_\fin$ is called \emph{unramified}
outside~$S$ if it is a finite linear combination of functions of the
form $g\cdot (1_{\MCR_T}\otimes f)$ for some smooth function~$f$
on~$\Adel_S$ and some $g\in\Adel\inv_\fin$.  We let $W_S\subseteq W$
be the subspace of functions that are unramified outside~$S$.  We
clearly have $W_S\subseteq W_{S'}$ for $S\subseteq S'$.  If we work
rationally, then a smooth function is unramified outside~$S$ if and
only if it is $\MCR\inv_T$\nbd{}invariant.  However, functions
in~$W_S$ must not only be $\MCR\inv_T$\nbd{}invariant, they also must
satisfy a stronger integrality condition.  Let
$$
L_m(y,S)\defeq \{g\in \GF\inv_\infty\mid gy-y\in x_m\MCR_S\}.
$$
Then $L_m(y)\subseteq L_m(y,S)$ for all~$S$.  It is straightforward to
see that~$W_S$ consists exactly of the $\MCR\inv_T$\nbd{}invariant
elements of~$W$ that satisfy $f(y)\in\bigcup_{m\in\N}
\Lambda(\dot{L}_m(y,S): \GF\inv_\infty/\Roots)$ for all
$y\in\MCR_\fin$.

For $a\in\N$, Let $W_a\subseteq W$ be the subspace generated by all
$W_S$ with $\# S\le a$.  This is an increasing filtration on~$W$.  We
use this filtration to compute the homology with coefficients in~$W$.
Let $W^{(0)}_a\subseteq W_a$ be the subspace that is spanned by
elements of the form $1_{\MCR_T} \otimes f_S$ for finite sets of
places~$S$ with $\# S=a$ and $\supp f_S\subseteq \MCR\inv_T$.  Let
$W^{(1)}\defeq W^{(0)}\cap W_{a-1}$.  Both $W^{(0)}_a$ and $W^{(1)}_a$
are representations of the compact open subgroup $\MCR\inv_\fin$ of
$\Adel\inv_\fin$.  Since this subgroup is open, compact induction from
$\MCR\inv_\fin$ to $\Adel\inv_\fin$ makes sense for representations on
Abelian groups.  Hence we can define a representation
$\cInd_{\MCR\inv_\fin}^{\Adel\inv_\fin} W^{(0)}_a/W^{(1)}_a$ of
$\Adel\inv_\fin$ and a canonical equivariant map
$$
\Phi_a\colon
\cInd_{\MCR\inv_\fin}^{\Adel\inv_\fin}
W^{(0)}_a/W^{(1)}_a \to W_a/W_{a-1}.
$$

\begin{lemma}  \label{lem:W_filtered_homology}
  The map~$\Phi_a$ is an isomorphism for all $a\in\N$.
\end{lemma}

\begin{proof}
  First we check surjectivity.  Let $f\in W_S$ for some set of finite
  places~$S$ with $\# S\le a$.  If $\# S<a$, then $f\in W_{a-1}$, so
  that we may assume $\# S=a$.  Write $f=g\cdot 1_{\MCR_T} \otimes
  f_S$.  We will show that $f\equiv f' \bmod W_{a-1}$ for a function
  $f'=g\cdot 1_{\MCR_T}\otimes f'_S$ with $\supp f'_S
  \subseteq\Adel\inv_S$.  Since $\Adel\inv_\fin$ acts transitively
  on~$\Adel\inv_S$, the function~$f'$ belongs to the range
  of~$\Phi_a$.  Hence~$\Phi_a$ must be surjective if we can modify~$f$
  in this fashion.

  If~$f$ is not yet supported in~$\Adel\inv_S$, there is $v\in S$ such
  that $f|_{x_v=0}$ does not vanish.  Write
  $\Adel_S=\GF_v\times\Adel_{S\setminus y}$ and define
  $$
  f_v(x_v,y) \defeq
  \begin{cases}
    f(0,y) & \text{if $\abs{x_v}_v\le1$;}
    \\
    0 & \text{if $\abs{x_v}_v>1$.}
  \end{cases}
  $$
  This is again an element of~$W$.  The reason for this is that the
  integrality condition at $(0,y)$ is stronger than the integrality
  condition at $(x_v,y)$ for $x_v\neq0$ because $L_m((0,y))\supseteq
  L_m((x_v,y))$.  Hence $f_v\in W_{S\setminus y}$.  Subtracting~$f_v$
  from~$f$, we obtain a function that vanishes whenever $x_v=0$.
  Proceeding in this fashion for all $v\in S$, we obtain a function
  supported in~$\Adel\inv_S$.  This establishes the surjectivity
  of~$\Phi_a$.

  Next we verify injectivity.  Let $A\subseteq \Adel_\fin$ be a closed
  and open subset.  Then multiplication by~$1_A$ is an operator on~$W$
  because the integrality condition is local.  If~$A$ is invariant
  under $\MCR\inv_\fin$, then this operator also preserves the
  subspaces~$W_S$ and hence the filtration.  It also acts on the
  source of~$\Phi_a$ in an evident fashion.  Both the source and
  target of~$\Phi_a$ carry a filtration defined by the condition that
  $\abs{x}_\fin \le n$ for all $x\in\supp f$.  Using multiplication
  operators, we obtain easily that~$\Phi_a$ is injective.
\end{proof}

The source of~$\Phi_a$ is a free module over $\GF\inv/\GF\inv_\infty$
because $\GF\inv/\GF\inv_\infty$ embeds in
$\Adel\inv_\fin/\MCR\inv_\fin$.  Hence it is relatively projective,
that is, projective for extensions that split as extensions of Abelian
groups.  The group $W^{(0)}/W^{(1)}$ contains torsion elements, so
that the source of~$\Phi_a$ is not free as an Abelian group.
Nevertheless, we can prove by induction on $a-b$ that $W_a/W_b$ is
relatively projective for all $b\le a$ because extensions of
relatively projective modules are again relatively projective.
Therefore, $W=\varinjlim W_a$ is relatively flat, so that
$H_m(\GF\inv/\GF\inv_\infty,W)=0$ for $m\neq0$.  This implies
$$
H_p(\GF\inv/\Roots,\K_\homogen)
\cong H_p(\GF\inv_\infty/\Roots,\K_\homogen)/\GF\inv
\cong W/\GF\inv.
$$
Our filtration $(W_a)$ also yields a filtration on $W/\GF\inv$, whose
subquotients are
$$
H_0(\GF\inv/\GF\inv_\infty, W_a/W_{a-1}) \cong
\CCINF(\Adel\inv_\fin/\GF\inv\MCR\inv_\fin) \otimes W_a^{(0)}/W^{(1)}_a.
$$
However, the extensions that come with this filtration are quite
complicated because of the torsion in the groups
$W_a^{(0)}/W^{(1)}_a$.

\end{document}